\documentclass{ws-ijm}

\begin{document}

\markboth{W. Wu, S. Zhang and Y.-Z. Zhang}
{Pointed Hopf algebras with classical  Weyl groups}

\catchline{}{}{}{}{}

\title{ON  NICHOLS (BRAIDED) LIE ALGEBRAS}

\author{WEICAI WU}

\address{Department  of Mathematics, Hunan University,
 Changsha  410082,  P.R. China }

\author{SHOUCHUAN ZHANG}

\address{Department  of Mathematics, Hunan University,
 Changsha  410082,  P.R. China }

\author{ YAO-ZHONG ZHANG}

\address{School of Mathematics and Physics, The University of Queensland,
 Brisbane 4072, Australia}

\maketitle

\begin{abstract}
We prove {\rm (i)} Nichols algebra $\mathfrak B(V)$ of vector space $V$ is
finite-dimensional if and only if Nichols braided Lie algebra $\mathfrak L(V)$ is finite-dimensional;
{\rm (ii)} If the rank of connected $V$ is $2$ and $\mathfrak B(V)$ is an arithmetic root system, 
then  $\mathfrak B(V) = F \oplus \mathfrak L(V);$ and {\rm (iii)}  if  $\Delta (\mathfrak B(V))$
is an arithmetic root system and there does not exist any $m$-infinity element with $p_{uu} \not= 1$
for any $u \in D(V)$, then  $\dim (\mathfrak B(V) ) = \infty$  if and only if there exists  $V'$,
which is twisting equivalent to $V$, such that $ \dim (\mathfrak L^ - (V')) = \infty.$
Furthermore we give an estimation of  dimensions of Nichols Lie algebras and two examples of Lie algebras 
which do not have maximal solvable ideals.
\end{abstract}



\newtheorem{Proposition}{Proposition}[section]
\newtheorem{Theorem}[Proposition]{Theorem}
\newtheorem{Definition}[Proposition]{Definition}
\newtheorem{Corollary}[Proposition]{Corollary}
\newtheorem{Lemma}[Proposition]{Lemma}
\newtheorem{Example}[Proposition]{Example}
\newtheorem{Remark}[Proposition]{Remark}

\maketitle 

\numberwithin{equation}{section}

\section{Introduction}\label {s0}
Nichols algebras play a prominent role in various areas of mathematics such as the theory of  pointed Hopf algebras and logarithmic quantum field theories.
Recently Heckenberger
\cite{He05} established a one-to-one correspondence between arithmetic root systems and Nichols algebras
of diagonal type as well as between twisted equivalence classes of arithmetic root systems
and generalized Dynkin diagrams. A great deal of attention has been paid to the question of finite-dimensionality of
Nichols algebras (see e.g. \cite {AHS08,AS10,An11,He05,He06a,He06b,WZZ14,WZZ15}).
The interest in this problem arose from the work of Andruskiewitsch and Schneider \cite{AS98}
on classification of finite dimensional (Gelfand-Kirillov) pointed Hopf
algebras which are generalizations of quantized  enveloping algebras of semi-simple Lie algebras.
In \cite{He05} Heckenberger classified braided vector spaces of diagonal type with
finite-dimensional Nichols algebras.

Let $\mathfrak B(V)$ be the Nichols algebra of vector space $V$.
Let  $\mathfrak L(V)$ ,   $\mathfrak L^ - (V)$ and $\mathfrak L_c(V)$ denote the  braided Lie algebras
 generated by $V$ in $\mathfrak B(V)$ under Lie operations $[x,    y] = yx  -  p_{yx}xy$,
$[x,    y]^ -  = xy  -  yx$ and $[x,    y]_c = xy  -  p_{xy}yx$,   respectively,
for any homogeneous elements
$x,    y \in \mathfrak B(V)$.   $(\mathfrak L(V),    [\ ])$,    $(\mathfrak L^ - (V),    [\ ]^ - )$ and
$(\mathfrak L_c(V),    [\ ]_c)$ are called Nichols braided Lie algebra,    Nichols Lie algebra
and Nichols braided m-Lie algebra of $V$,    respectively.
It is clear that $(\mathfrak L(V),    [ \ ])$ and $(\mathfrak L_c(V),    [ \ ]_c)$
are equivalent as vector spaces. If $\mathfrak B(V)$ is finite dimensional
then $\mathfrak B(V)$ is nilpotent,    so
$(\mathfrak L(V),    [\ ])$ and $(\mathfrak L^ - (V),    [\ ]^ - )$ also are nilpotent.

In \cite{WZZ15}, we studied the relationship between Nichols braided Lie algebras
and  Nichols algebras and provided a new method to determine when a Nichols algebra is finite dimensional.
It was shown there that Nichols algebra $\mathfrak B(V)$ is
finite-dimensional if and only if Nichols braided Lie algebra $\mathfrak L(V)$ is finite-dimensional
if there does not exist any $m$-infinity element in $\mathfrak B(V)$, and that
$\mathfrak B(V) = F\oplus \mathfrak L(V)$ does not hold if the rank of connected $V$ is higher than  $3$
and $\mathfrak B(V)$ is an arithmetic root system.

In this paper we show that the condition ``there does not exist any $m$-infinity element in
$\mathfrak B(V)$" can be dropped and Nichols algebra $\mathfrak B(V)$ is
finite-dimensional if and only if Nichols braided Lie algebra $\mathfrak L(V)$ is finite-dimensional.
Furthermore we prove that $\mathfrak B(V) = F \oplus \mathfrak L(V)$ if the rank of connected $V$ is equal to $2$
and $\mathfrak B(V)$ is an arithmetic root system, and that $\dim (\mathfrak B(V) ) = \infty$
if and only if there exists  $V'$,  which is twisting equivalent to $V$, such that
$ \dim (\mathfrak L^ - (V')) = \infty$ if  $\Delta (\mathfrak B(V))$ is an arithmetic root system
and there does not exist any $m$-infinity element with $p_{uu} \not= 1$ for any $u \in D(V)$.
Finally we give an estimation of  dimensions of Nichols Lie algebras.

This paper is organized as follows. In section \ref{s2} we prove that Nichols algebra $\mathfrak B(V)$ is
finite-dimensional iff Nichols braided Lie algebra $\mathfrak L(V)$ is finite-dimensional.
In section \ref {s3} we show that a generalized Dynkin diagram of $V$ with rank $n$ is a   complete diagram with $q_{ii} \not=1$ for any $1\le i \le n$ iff  $\mathfrak B(V) = F \oplus \mathfrak L(V)$. This implies that
if the rank of connected $V$ is equal to $2$ and $\mathfrak B(V)$ is an
arithmetic root system,   then  $\mathfrak B(V) = F \oplus \mathfrak L(V).$
  In  section \ref{s5} we show
that if $\Delta (\mathfrak B(V))$ is an arithmetic root system
and there does not exist any $m$-infinity element with $p_{uu} \not= 1$ for any $u \in D(V)$,  then $\mathfrak B(V)$ is
finite-dimensional iff $\mathfrak L^ - (V)$ is finite-dimensional.
We also give an estimation of dimensions of Nichols Lie algebras.
In section \ref{solvable ideals} we find two examples of Lie algebras which have no maximal solvable ideals. In  the appendix we prove that there  does not exist any $m$-infinity element in Nichols algebra
 $\mathfrak B(V) $ with arithmetic root system $\Delta (\mathfrak B(V))$  over finite cyclic groups.

Throughout,   $\mathbb Z =: \{x | x \hbox { is an  integer }\},$
$\mathbb N_0 =: \{x | x \in \mathbb Z,    x\ge 0\},$
$\mathbb N =: \{x | x \in \mathbb Z,    x>0\}.$  $F$ denotes the base field of characteristic zero. $ N_{k} := {\rm ord } (p_{k  k})$ when $p_{k k} \not=1$; $ N_{k} := \infty$  when $p_{k k} =1$. $D  = D(V)=: \{[u] \mid [u] \hbox { is a hard super-letter }\}$.
If $[u] \in D$ and
${\rm ord }  (p_{u u}) =m>1$ with $h_u = \infty$, then $[u]$ is called an $m$-infinity element.
Other notations are the same as in  \cite {AHS08} and \cite {WZZ15}. Throughout, braided vector space $V$ is connected and  of diagonal type with basis $x_1, x_2, \cdots, x_n$ and $C(x_i \otimes x_j) = q_{ij} x_j \otimes, x_i$, $n>1 $.

\section{Relationship between Nichols algebras and Nichols Lie algebras}\label {s2}
In this section we show that Nichols algebra $\mathfrak B(V)$ is
finite-dimensional if and only if Nichols braided Lie algebra $\mathfrak L(V)$ is finite-dimensional.

\begin {Lemma} \label {1.1}  If  homogeneous element   $u\in \mathfrak{L}(V)$ and  ${\rm ord }  (p_{uu})>2$,   then $u^{m}\in\mathfrak{L}(V),  \forall m\in\mathbb N$.
\end {Lemma}

\noindent {\bf Proof.} We prove this by induction on $m$. If $m=1$,   the claim holds obviously. Use the induction hypothesis,   we have  $u,  u^{2},  \cdots,  u^{m}\in\mathfrak{L}(V)$. If $u^{m +  1}\notin\mathfrak{L}(V)$,   
by Lemma 4.12 in \cite {WZZ15} and
$\left \{
\begin
{array} {llll}
\ [u,  u^{m}]=(1 - p_{uu}^{m})u^{m +  1}   \ \ \ \ \ \ \ \\
\ [u^{2},  u^{m - 1}]=(1 - p_{u  u}^{2(m - 1)})u^{m +  1}  \ \ \ \ \ \ \
\end {array} \right. , $
we know $1 - p_{u u}^{m}=0,  1 - p_{u  u}^{2(m - 1)} =0$,   i.e.   $p_{u  u}^{2}=1$,   which is a contradiction to ${\rm ord }  (p_{uu})>2$. \hfill $\Box$

 \begin {Theorem} \label {1.2} If $\mathfrak B(V) $ is connected Nichols algebra of diagonal type with $\dim V>1$,   then $\mathfrak B(V)$ is
finite-dimensional if and only if $\mathfrak L(V)$ is finite-dimensional.
\end {Theorem}

\noindent {\bf Proof.} If there does not exist any $m$-infinity element in $\mathfrak B(V)$,   it follows from Theorem 4.11 in \cite {WZZ15}. We now show that the theorem also holds if there exists an $m$-infinity element $u$ in $\mathfrak B(V)$. In this case $\mathfrak B(V)$ is
infinite-dimensional.
 If ${\rm ord }  (p_{uu})=m>2$,  then $u^{k}\in\mathfrak{L}(V)$ by Lemma \ref {2.1},  for $\forall k\in\mathbb N$. Consequently, $\mathfrak L(V)$ is infinite-dimensional. Assume  that $p_{uu}^{2}=1$ and  $\mathfrak L(V)$ is finite-dimensional, then   we know $\exists v\in D$ such that $p_{uv}p_{vu}\neq1$ with $v\neq u$ by  Proposition 4.10 in \cite{WZZ15}.  We prove $u^{k}v\in\mathfrak{L}(V)$ by induction on $k$. If $k=1$,   the claim holds by Lemma 4.12 in\cite  {WZZ15}. Using the induction hypothesis we have $u^{k}v\in\mathfrak{L}(V)$. If $u^{k +  1}v\notin\mathfrak{L}(V)$,   then $p_{uu}^{2k}p_{uv}p_{vu}=1$ by Lemma 4.12 in \cite  {WZZ15}. which is a contradiction to $p_{uu}^{2}=1$ and $p_{uv}p_{vu}\neq1$.  This has shown $u^{k }v\in\mathfrak{L}(V)$ for any $k \in \mathbb N.$  It is clear that  $\{u^k v \mid k\in \mathbb N\}$ or $\{v u^k \mid k\in \mathbb N\}$ is a subset of restricted  PBW basis. Consequently, $\mathfrak L(V)$ is infinite-dimensional, which is a contradiction.
\hfill $\Box$

\section{Conditions for $\mathfrak B(V)=F\oplus \mathfrak L(V).$}\label {s3}
In this section we prove that a generalized Dynkin diagram  of $V$ with rank $n$ is a   complete diagram with $q_{ii} \not=1$ for any $1\le i \le n$ if and only if  $\mathfrak B(V) = F \oplus \mathfrak L(V)$. This implies that
if the rank of connected $V$ is equal to $2$ and $\mathfrak B(V)$ is an
arithmetic root system,   then  $\mathfrak B(V) = F \oplus \mathfrak L(V).$

\begin{Lemma}\label{2.1} Set $\widetilde{p}_{uv}:=p_{uv}p_{vu}$.  Let $a,  b,  c,  d,  e,  f$ denote  $1 - \widetilde {p}_{vw},  1 -\widetilde { p}_{uw},  1 - \widetilde{p}_{vu},  1 - \widetilde {p}_{vu}\widetilde{p}_{uw},  1 - \widetilde{p}_{vu}\widetilde{p}_{vw},  1 - \widetilde{p}_{vw}\widetilde{p}_{uw}$,   respectively.
If  homogeneous elements $u,  v,  w\in\mathfrak{L}(V)$,   then $uvw,  uwv,  vwu,  vuw,  wuv,  wvu\in\mathfrak{L}(V)$ if one of the following is fulfilled:

{\rm (i)} $a\neq0,  b\neq0$;

{\rm (ii)} $a\neq0,  c\neq0$;

{\rm (iii)} $a\neq0,  d\neq0$;

{\rm (iv)} $b\neq0,  c\neq0$;

{\rm (v) } $b\neq0,  e\neq0$;

{\rm (vi)} $c\neq0,  f\neq0$.
\end {Lemma}

\noindent {\bf Proof.}
 Analogously to a formula of Kharchenko on braided commutators in \cite {Kh99b}, the following equation holds: $[[u,v],w] = [u,[v,w]] + p_{wv}[[u,w],v] + p_{vu}(p_{wv}p_{vw}-1)[u,w]v$.

{\rm (i)} Assume now that $b\neq0$ and $a\neq0$. Then $[u,w]v\in \mathfrak L(V)$ and hence $v[u,w]\in \mathfrak L(V)$ by this formula since $a\neq0$. If $c\neq0$, then $[w,u]v\in \mathfrak L(V)$ and $v[w,u]\in \mathfrak L(V)$.
 there exist $\alpha, \beta, \alpha ' , \beta ' \in F$ such that $wu = \alpha [u, w] + \beta [w, u]$ and  $uw = \alpha' [u, w] + \beta' [w, u]$ since $b\not=1$. Consequently, $uw v = \alpha [u, w]v + \beta [w, u]v \in \mathfrak L(V)$ and  $wuv = \alpha' [u, w]v + \beta' [w, u]v\in \mathfrak L(V).$ Similarly, $uvw,    vuw,  vwu, wvu \in\mathfrak{L}(V)$.
 If $c=0$, then $d\neq0$ and $e\neq0$, then $(uw)v,(wu)v,v(uw),v(wu), (wv)u,u(vw)\in \mathfrak L(V)$ by 
 Lemma 4.12 in \cite {WZZ15}.

  Similarly, we can prove {\rm (ii)} and {\rm (iv)}.

  {\rm (iii)}
Assume now that $a\neq0$ and $d\neq0$. Then $b\neq0$ or $c\neq0$. In the first case {\rm (iii)} follows from {\rm(i)}. In the second, {\rm (iii)} also follows from {\rm (ii) }.

 Similarly, we can prove {\rm (v)} and {\rm (vi)}.
 \hfill $\Box$

\begin{Lemma}\label{2.2} Assume that homogeneous elements $u$ and $v$ are  in $ {\mathfrak L}(V)$.

{\rm (i)} Then $uv\in\mathfrak{L}(V)\Longleftrightarrow vu\in\mathfrak{L}(V)$.

{\rm (ii)} If $  p_{uv}p_{vu}\neq1$ and $p_{uu}^{2}=1$,   then $u^{m}v,  u^{m - 1}vu,  \cdots,  vu^{m}\in\mathfrak{L}(V),  \forall\ m\in\mathbb{N}$.

{\rm (iii)} Assume that $p_{uu}^{2}\neq1$ and that $u, v, uv\in\mathfrak{L}(V)$. Then $u^{k}vu^{l}\in\mathfrak{L}(V)$ for any $k,l\in \mathbb N_0$.

{\rm (iv)}  If $p_{uv}p_{vu}\neq1$ or $p_{u u}^2 \not=1$
with $uv\in \mathfrak{L}(V)$,   then $u^{m}v,  u^{m - 1}vu,  \cdots,$
$vu^{m} \in \mathfrak{L}(V),  \forall\ m\in \mathbb{N}$.
\end {Lemma}

\noindent {\bf Proof.} {\rm (i)} It is clear.

{\rm (ii)} We prove this by induction on $m$. set $u:=u,  v:=u$ and let $w:=u^{m - 2}v,  u^{m - 3}vu,  \ldots,  $

\noindent $vu^{m - 2},  m\geq2$,   respectively. then $a=1 - p_{uu}^{2(m - 2)}p_{uv}p_{vu}\neq0,  b=1 - p_{uu}^{2(m - 2)}p_{uv}p_{vu}\neq0$,   then $u^{m}v,  u^{m - 1}vu,  \cdots,  vu^{m}\in\mathfrak{L}(V)$ by Lemma \ref {2.1}.

{\rm (iii)} The proof can be done by induction on $k+l$. (i) and the induction hypothesis imply that it suffices to prove that $u^{m}v \in\mathfrak{L}(V)$ for any $m\in \mathbb N$. This again is done by induction on $m$. The claim holds by assumption for $m=1$.

Let $m\geq1$ such that $u^{k}v \in\mathfrak{L}(V)$ for any $k\in \{0,1,\cdots,m\}$. Since $p_{uu}^{2}\neq1$, Lemma \ref {1.1} implies that $u^{k} \in\mathfrak{L}(V)$ for any $k\geq1$. If $p_{u^{m}v,u}p_{u,u^{m}v}\neq1$, then $u^{m+1}v \in\mathfrak{L}(V)$.  By the same reason, if $p_{u^{m-1}v,u^{2}}p_{u^{2},u^{m-1}v}\neq1$, then $u^{m+1}v \in\mathfrak{L}(V)$. If both inequalities fail, then $p_{uv}p_{vu}=p_{uu}^{-2m}$ and $p_{uu}^{4}=1$. Since $p_{uu}^{2}\neq1$, we conclude that $p_{uu}^{2}=-1$ and $p_{uv}p_{vu}=(-1)^{-m}$.

Assume now that $m$ is odd. Apply Lemma \ref {2.1} to the triple $(u,u^{m},v)$. Since $p_{u,u^{m}}p_{u^{m},u}=p_{uu}^{2m}=-1$ and $p_{uv}p_{vu}=-1$, Lemma \ref {2.1} implies that
$u^{m+1}v \in\mathfrak{L}(V)$.

Assume now that $m$ is even. Then $p_{uv}p_{vu}=1$ and $m\geq2$. Apply Lemma \ref {2.1} to the triple $(u,u^{m-1},uv)$. Since $p_{u,u^{m-1}}p_{u^{m-1},u}=p_{uu}^{2m-2}=-1$ and $p_{u,uv}p_{uv,u}=p_{uu}^{2}p_{uv}p_{vu}=-1$, Lemma \ref {2.1} implies that $u^{m+1}v \in\mathfrak{L}(V)$. This proves the claim.

{\rm (iv)} It is clear by {\rm (ii)} and {\rm (iii)}. \hfill $\Box$

\begin {Lemma} \label {2.3} Let $m\in\mathbb N$ and  for any $1 \le i \le m$ let $w_i$ be a homogeneous element in $ {\mathfrak L}(V)$ with   $w_i ^{2 } =0$ when $p_{w_i, w_i} ^2=1$. If $p _{w_i,   w_{j}} p _{w_{j},   w_i} \not=1$   when  $w_i \not= w_j$ for any $1\le i\not= j \le m$,
    then $W: =  \prod \limits_{i=1} ^m w_i \in \mathfrak L(V) \oplus F   $.

\end {Lemma}
\noindent {\bf Proof.}
  It is clear for $m=1.$ For $m=2,$ if $w_1= w_2$ and $p_{w_1, w_2} p_{w_2, w_1} =1$, then $p_{w_1, w_1} ^2 =1$ and $w_1^2 = w_1w_2 =0 \in \mathfrak L(V) \oplus F .$
  If  $p_{w_1, w_2} p_{w_2, w_1} \not=1$, then $ w_1w_2 \in \mathfrak L(V) \oplus F$ by Lemma 4.12 in \cite  {WZZ15}.

Assume  $m>2.$  If the claim does not hold,   then there exists a minimal $t$ such that $W :=  \prod \limits_{i=1} ^t w_i \notin \mathfrak L(V). $ Obviously, $t>2$, $W\not=0$ and $w_1w_2\cdots w_{t-2}, w_1w_2\cdots w_{t-1} \in \mathfrak L(V).$ Let $u:= w_{t-1},$ $ v:= w_t, $ $w : = w_1w_2\cdots w_{t-2}$. If $u=v$ and $p_{uv}p_{vu} =1,$ then $p_{uu}^2 =1$, which implies $uu= uv =0.$ This contradicts  to $W :=  \prod \limits_{i=1} ^t w_i \notin \mathfrak L(V) .$ Consequently, $p_{uv} p_{vu} \not=1.$
It follows from Lemma 4.12 in \cite  {WZZ15} that $\widetilde{p} _{v, wu} =1$ .
By  Lemma \ref {2.1}, we have that 
$\widetilde{p} _{v, w} =1$, which  implies $\widetilde{p} _{uv} =1$. This is a contradiction.
 Here $\widetilde{p} _{uw} := p_{uw} p_{wv}.$
\hfill $\Box$

\begin {Lemma} \label {2.4} Assume that $m \in \mathbb N$ and $w_i$ is a homogeneous element in $ {\mathfrak L}(V)$ with   $w_i ^{2 } =0$ when $p_{w_i, w_i} ^2=1$ for any $1 \le i \le m$.
If $p _{w_i,   w_{j}} p _{w_{j},   w_i} \not=1$   when  $w_i \not= w_j$ for any $1\le i\not= j \le m$,    then $W :=  \prod \limits_{i=1} ^m w_i^{a_i} \in \mathfrak L(V) \oplus F    $ for any  $a_i \in \mathbb N_0$.

\end {Lemma}
\noindent {\bf Proof.} We show this by following several steps.

{\rm (i)}
Set $N(W) := \sum _{i = 1}^m a _i $. We show $W \in \mathfrak L(V)\oplus F$  by induction on $N(W).$  It is clear when $W=0$. We assume $W \not=0$ from now on.  When $N(W) =0,  $ $W=1 \in \mathfrak L(V)\oplus F$. When $N(W) =1,  $ $W \in \mathfrak L(V)\oplus F$.
Now assume $N(W) >1.$

{\rm (ii)} If $m=1$, $W=w_1^{a_1} \in \mathfrak L(V)\oplus F$ by Lemma \ref {1.1}.

{\rm (iii)} If $ a_i =1$ for $1\le i \le m, $ then $W \in \mathfrak L(V)\oplus F$ by Lemma \ref {2.3}.

{\rm (iv)}
 If  there exists $1 \le i_0 \le  m$ such that $a_ {i_0}>1$,    let $z_1 = w_{i_0},   z_2 =  w_{i_0 +  1} ^{a_{i_0 +  1}} \cdots w_{m}^{a_{m}} w_1 ^ {a_1} \cdots w_{i_0 - 1} ^{a_{i_0 - 1}}.$
By induction hypothesis and Lemma \ref {2.2} {\rm (i)},   $z_1,   z_2,   z_1z_2 \in \mathfrak L(V)$.  By Lemma \ref {2.2}(iv),   $ W'= z_1 ^{a_{i_0}} z_2\in {\mathfrak L} (V)$. Using induction again, we have $ w_{i_0}^{a_{i_0}} \cdots w_m ^{a_m}, w_1 ^{a_1} \cdots w _{i_0-1}^{a_{i_0-1}} \in {\mathfrak L} (V)$. By Lemma \ref {2.2} {\rm (i)}  we have that $W
 = w_1 ^{a_1} \cdots w _{i_0-1}^{a_{i_0-1}} w_{i_0}^{a_{i_0}} \cdots w_m ^{a_m} \in {\mathfrak L} (V)$.
\hfill $\Box$

\vskip.1in
A generalized Dynkin diagram is called  a complete diagram if  there exists an edge between  any two vertexes.

\begin {Theorem} \label {2.7} A generalized Dynkin diagram  of $V$ with rank $n$ is a complete diagram with $q_{ii} \not=1$ for any $1\le i \le n$ if and only if  $\mathfrak B(V) = F \oplus \mathfrak L(V).$

\end {Theorem}
\noindent {\bf Proof.}  The necessity. Let  $w_i \in \{ x_1, x_2, \cdots, x_n\}$   for any $1\le i\le m$. Then $p_{w_i, w_{j}} p_{w_{j}, w_{i}} \not= 1 $  when $w_i \not= w_j$ for any $1\le i \not= j \le m$
and
$\prod \limits_{i=1} ^m w_i^{a_i} \in {\mathfrak L } (V) \oplus F$ for any $a_i \in \mathbb N_0$ by Lemma \ref {2.4}. Consequently, $\mathfrak B(V) = F \oplus \mathfrak L(V).$

The sufficiency. Obviously $ q_{ii} \not=1$ since $x_i^2 \in  F \oplus \mathfrak L(V)$
 for $1\le i\le n.$ By Lemma 5.2 in \cite {WZZ15}, we have that the  generalized Dynkin diagram  of $V$ is a   complete diagram.

 \hfill $\Box$

We immediately have

\begin {Corollary} \label {2.9} Assume that $\Delta (\mathfrak B(V))$ is an arithmetic root system and $V$ is connected. Then $\mathfrak B(V) = F \oplus \mathfrak L(V)$ if the rank of  $V$ is $2$.


\end {Corollary}

\begin {Lemma} \label {2.5} Assume that $w_i$ is a homogeneous element in $ {\mathfrak L}(V)$ with    $w_i ^{2 } =0$ when $p_{w_i, w_i}^2 =1$ for any $1 \le i \le m$.
 If $\prod \limits_{j=1} ^r w_{i_j} \in \mathfrak L(V) \oplus F $ for any distinct  $1\le i_2, i_2, \cdots, i_r \le m$,
    then $W :=  \prod \limits_{i=1} ^m w_i^{a_i} \in \mathfrak L(V) \oplus F   $ for any  $a_i \in \mathbb N_0$.

\end {Lemma}
\noindent {\bf Proof.}
Set $N(W) := \sum _{i = 1}^m a _i $. We show $W \in \mathfrak L(V)\oplus F$  by induction on $N(W).$  It is clear when $W=0$. We assume $W \not=0$ from now on.  When $N(W) =0,  $ $W=1 \in \mathfrak L(V)\oplus F$. When $N(W) =1,  $ $W \in \mathfrak L(V)\oplus F$.
Now assume $N(W) >1.$

{\rm (ii)} If $m=1$, $W=w_1^{a_1} \in \mathfrak L(V)\oplus F$ by Lemma \ref {1.1}.

{\rm (iii)} If $ a_i = 1$ for $1\le i \le m, $ then $W \in \mathfrak L(V)\oplus F$ by assumption.

{\rm (iv)}
 If  there exists $1 \le i_0 \le  m$ such that $a_ {i_0}>1$,    let $z_1 = w_{i_0},   z_2 =  w_{i_0 +  1} ^{a_{i_0 +  1}} \cdots w_{m}^{a_{m}} w_1 ^ {a_1} \cdots w_{i_0 - 1} ^{a_{i_0 - 1}}.$
By induction hypothesis and Lemma \ref {2.2} {\rm (i)},   $z_1,   z_2,   z_1z_2 \in \mathfrak L(V)$.  By Lemma \ref {2.2}(iv),   $  z_1 ^{a_{i_0}} z_2\in {\mathfrak L} (V)$. Using induction again, we have $ w_{i_0}^{a_{i_0}} \cdots w_m ^{a_m}, w_1 ^{a_1} \cdots w _{i_0-1}^{a_{i_0-1}} \in {\mathfrak L} (V)$. By Lemma \ref {2.2} {\rm (i)}  we have that $W
 = w_1 ^{a_1} \cdots w _{i_0-1}^{a_{i_0-1}} w_{i_0}^{a_{i_0}} \cdots w_m ^{a_m} \in {\mathfrak L} (V)$.
\hfill $\Box$

\begin {Lemma} \label {2.6}
 If homogeneous element $u$ in $ {\mathfrak L}(V)$   and $p_{u,   u} \not=1,  $ then $u^2\in {\mathfrak L} (V). $

\end {Lemma}
\noindent {\bf Proof.}  $[u,   u] = u^2  -  p_{u,   u} u^2 = (1 -  p_{u,   u}) u^2$,   which implies $u^2\in {\mathfrak L} (V)$. \hfill $\Box$

\begin {Proposition} \label {2.8} Assume that  $u$ and $v$ are  homogeneous elements in $ {\mathfrak L}(V)$ and  $(p_{uv,   uv})^2 \not=1$  with $  uv, uvu, vuv \in \mathfrak L(V)$.
  If $ u^2 =0  $ when $p_{uu}^2 =1  $ and $v^2 =0 $ when $p_{vv}^2 =1$,      then $W =  \prod _{i=1} ^mu^{a_i} v^{b_i} \in \mathfrak L(V)   $ for any $a_i,   b_i \in \mathbb N_0,  $ $m\in \mathbb N.$

\end {Proposition}
\noindent {\bf Proof.}
  Set $N(W) := \sum _{i = 1}^m (a _i  +  b_i)$. We show $W \in \mathfrak L(V)\oplus F$  by induction on $N(W).$ When $N(W) =0,  $ $W=1 \in \mathfrak L(V)\oplus F$. When $N(W) =1,  $ $W \in \mathfrak L(V)\oplus F$. Assume $N(W) >1.$ Set
$c_{2i - 1} = a_i,   c_{2i}= b_i$,   $y_{2i - 1} = u,   y_{2i} =v $ for $1\le i \le m.$

{\rm (i)} $(uv) ^m u,$ $ (vu) ^m v,$ $ (vu) ^m , (uv) ^m \in \mathfrak L(V)$  for any $m \in \mathbb N.$ In fact, $ (vu) ^m, (uv) ^m \in  \mathfrak L(V)$ for any $m \in \mathbb N$ by Lemma \ref {1.1}. $(uv) ^m u,$ $ (vu) ^m v  \in \mathfrak L(V)$ by Lemma \ref {2.2} {\rm (iv)}.

{\rm (ii)} By Lemma \ref {1.1},   $y_i ^{c_i} \in  \mathfrak L(V) $  for $1\le i \le m$. Now we assume that there are two non-zero elements in $\{c_i \mid 1\le i \le m\}.$

{\rm (iii)} If  there exists $1 \le i_0 \le  2m$ such that $c_ {i_0}>1$,    let $z_1 = y_{i_0},   z_2 =  y_{i_0 +  1} ^{c_{i_0 +  1}} \cdots y_{2m}^{c_{2m}} y_1 ^ {c_1} \cdots y_{i_0 - 1} ^{c_{i_0 - 1}}.$
By induction,   $z_1,   z_2,   z_1z_2 \in \mathfrak L(V)$.  By Lemma \ref {2.2}(iv),   $ W'= z_1 ^{c_i} z_2\in {\mathfrak L} (V)$. Using induction again,   we have $W \in {\mathfrak L} (V)$.
\hfill $\Box$

\begin {Proposition} \label {2.10}  $V$ is a quantum linear space if and only if  ${\rm span } \{ x_i^{a_i} \mid  0 < a_i < N_i, 1\le i \le n\} = \mathfrak L(V)$. \end {Proposition}

\noindent {\bf Proof.} The sufficiency. It is clear $[x_i, x_j] =0 $ for any $i \not= j, $  i.e. $p_{ij}p_{ji} =1$ for any $i \not= j. $ Consequently, $V$ is a quantum linear space.
The necessity. It is clear $x_i ^{a_i} \in \mathfrak L(V). $  It is enough to show that
$b(x_{i_1}, x_{i_2}, \cdots, x_{i_r})=0$  when there exist two elements in $x_{i_1}, \cdots, x_{i_r}$ are different, where $b(x_{i_1}, x_{i_2}, \cdots, x_{i_r})$ denote a method of adding  bracket on  $x_{i_1}, x_{i_2}, \cdots, x_{i_r}$. We show this by induction on $r$. It is clear when $r=2$. When $r>2$, let  $b(x_{i_1}, x_{i_2}, \cdots, x_{i_r}) = [u_1, u_2]$ with $u_1 = b(x_{i_1}, \cdots, x_{i_s} )$. If   $u_1=0$ or
 $u_2 =0$, then $b(x_{i_1}, x_{i_2}, \cdots, x_{i_r}) =0$. If $u_1\not=0$ and $u_2\not=0$, then  $x_{i_1} = \cdots = x_{i_s}$ and $x_{i_{s+1}} = \cdots = x_{i_r}$ by induction. It follows  $b(x_{i_1}, x_{i_2}, \cdots, x_{i_r})=0$. from  braided Jacobi identity. \hfill $\Box$

\section{ Nichols Lie algebra   $\mathfrak L^ - (V)$}\label {s5}
In this section it is proved that if $\Delta (\mathfrak B(V))$ is an arithmetic root system
and there does not exist any $m$-infinity element with $p_{u,u} \not= 1$ for any $u \in D(V)$ then
  $\dim (\mathfrak B(V) ) = \infty$  if and only if there exists  $V'$,  which is twisting equivalent to $V$, such that $ \dim (\mathfrak L^ - (V')) = \infty.$
Furthermore, we give an estimation of  dimensions of Nichols Lie algebras.

If  $A$ is an associative algebra and we define  $[a, b]^- = ab -ba$ for any $a, b \in A$, then $(A, [\  ]^-)$ is a Lie algebra, which is written as $A^-.$
Recall that Nichols Lie algebra $\mathfrak L^ - (V)$  is the Lie algebra generated by $V$ in $\mathfrak B (V)^-$ (see \cite {WZZ15}).

For braided vector space $V$, we can define a directed diagram as follows. We   use  an arrow from $i$ to $j$ to denote $p_{ij} \not=1$:
 $$\begin{picture}(250,    15) \put(27,   1){\makebox(0,
0)[t]{$\bullet $}} \put(60,    1){\makebox(0,   0)[t]{$\bullet$}}
\put(28,    - 1){\line(1,    0){33}}\put(20,
 - 20){i} \put(58,   - 20){j}\put(56,     - 4.5){$\rangle$}
\end{picture} $$
We use  a  bi-arrow between  $i$ and  $j$ to denote
$p_{ij} \not=1$ and $p_{ji} \not=1$:
$$
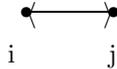
\begin{picture}(250,    15) \put(27,   1){\makebox(0,
0)[t]{$\bullet$}} \put(60,    1){\makebox(0,   0)[t]{$\bullet$}}
\put(28,    - 1){\line(1,    0){33}}\put(20,
 - 20){i} \put(58,   - 20){j}\put(56,     - 4.5){$\rangle$}\put(27,    - 4.5){$\langle$}
\end{picture}$$  \ \ \\
This directed diagram   is   called  a directed generalized Dynkin diagram. Furthermore, we add a line between $i$ and $j$  when    $p_{ij}p_{ji} \not= 1$. In this case, we call this diagram a mixed  generalized Dynkin diagram.
 We can study Nichols Lie algebras by means of the directed generalized Dynkin diagrams.

Let
$l_{u}^{0}[v]=:[v],  l_{u}^{i}[v]=:[[u],  l_{u}^{i - 1}[v]],  r_{u}^{0}[v]:=[v],  r_{u}^{i}[v]:=[r_{u}^{i - 1}[v],  [u]],  i\geq1$
for any $u,   v\in \mathfrak B(V)$. similarly,  Let
$\overline{l}_{u}^{0}[v]^ - =:[v]^ - ,  \overline{l}_{u}^{i}[v]^ - =:[[u]^ - ,  \overline{l}_{u}^{i - 1}[v]^ - ]^ - ,  \overline{r}_{u}^{0}[v]^ - :=[v]^ - ,
\overline{r}_{u}^{i}[v]^ - $ $:=[\overline{r}_{u}^{i - 1}[v]^ - ,  $ $[u]^ - ]^ - ,  i\geq 1$
for any $u,   v\in \mathfrak B(V)$. Let $(k)_a : = 1  +  a  +  a^2  +   \cdots  +   a^{k - 1};$ $(k)_a! := (1)_a (2)_a \cdots (k)_a.$

Recall the dual $\mathfrak B(V^*) $ of Nichols algebra $\mathfrak B(V)$ of rank $n$ in Section 1.3 of \cite {He05} and \cite{He06b}. Let $y_{1},   y_2,   \cdots,   y_n$ be a
dual basis of $x_{1},   x_2,   \cdots,   x_n$. $\delta(y_i)=g_i ^{ - 1}\otimes y_i$,   $g_i\cdot y_j=p_{ij}^{ - 1}y_j $ and $\Delta (y_i)=g_i^{ - 1}\otimes y_i +  y_i\otimes1$. There
exists a bilinear map $<,   >$ from $(\mathfrak B(V^*)\# kG)\times\mathfrak B(V)$ to $\mathfrak B(V)$ such
that $<y_i,  uv> = <y_i,  u>v  +  g_i^{ - 1}.u<y_i,      v>$ and $<y_i,  <y_j,  u>> = <y_iy_j,  u>$  for any $u,  v\in\mathfrak B(V)$.
Furthermore,  for any $u\in\oplus _{i=1}^\infty\mathfrak B(V)_{(i)}$,   one has that $u=0$ if and only if $<y_i,  u>=0$ for
any $1\leq i\leq n.$ Let $i$ denote $x_{i}$ in short,   sometimes.

\begin {Lemma} \label {3.1} Let  $a:= p_{ii}^{ - 1}$,   $ b:= p_{ij}^{ - 1}$,   $c:= p_{ji}^{ - 1}$. Then

{\rm (i)} $<y_j,   [x_i,   x_j]^ - > = (p_{ji} ^{ - 1} - 1) x_i,   $ $<y_i,   [x_i,   x_j]^ - > = (1 - p_{ij} ^{ - 1}) x_j,   $

{\rm (ii)} (See Lemma 1.3.3 in \cite{He05}) $<y_i,   x_i^m> = (m)_a x_i^{m - 1}. $

{\rm (iii)} $\overline{l}_{i}^{m}[j]^ -  = \sum \limits _{k=0}^{m} ( - 1)^k \left(       \begin{array}{ccccccccccccccccc} m\\
k
                   \end{array}\right) x_i ^{m - k}x_j x_i ^k.$

{\rm (iv)}  $ <y_{i},  \overline{l}_{i}^{m}[j]^ - >=\sum \limits
_{k=0}^{m - 1}( - 1)^{k}\{ \left (\begin{array}{cc}
m\\
k\ \end{array} \right  )(m - k)_a -  \left (\begin{array}{cc}
m\\
k +  1\
\end{array} \right )((k +  1)_a)a^{m - k - 1}b\}x_{i}^{m - k - 1}x_{j}x_{i}^{k}$.

{\rm (v)} $<y_{i}^{m}y_{j},  \overline{l}_{i}^{m}[j]^ - >=   (c - 1)^m (m)_a!$.

{\rm (vi)}  $<y_{i}^{m - 1}y_{j}y_{i},  \overline{l}_{i}^{m}[j]^ - >
$ $=  \frac {(m - 1)_a!} {c(1 - a)} ( (c - 1)^m (1 - ba^m c)  +
 (ac - 1)^m (cb  - 1))$,   when $a\not=1$

 {\rm (vii)} $\overline{l}_{i}^{m}[j]^ -  \not= 0,  $  when  $ c\not= 1$ and $ {\rm ord } (a) > m.$

{\rm (viii)} $\overline{l}_{i}^{m}[j]^ - \neq0$,   when  $l_{i}^{m}[j]\neq0$ and $a \not= 1$.

{\rm (ix)}  $x_i^m \notin \mathfrak L^ -  (V) $ for $1<m$ and $x_i^m \not= 0.$

{\rm (x)} $\dim (\mathfrak L^ - (V)) = \infty$ when there exists $i$ and $j$ with  $i \not= j$,   $p_{ji} \not= 1$ and ${\rm ord } (p_{ii}) = 1$ or $\infty$.
\end {Lemma}

\noindent {\bf Proof.}
{\rm (i)} It can be proved by simple computation.

{\rm (iii)} and {\rm (iv)} They  can be proved by induction.

{\rm (v)} It follows from Part {(ii)} and Part {(iii)}.

{\rm (vi)}  \begin {eqnarray*} && <y_{i}^{m - 1}y_{j}y_{i},  \overline{l}_{i}^{m}[j]^ - >\\
&=&
\sum \limits _{k=0}^{m - 1}( - 1)^{k}\{(\begin{array}{cc}
m\\
k\ \end{array})(a^{m - k} - 1) - (\begin{array}{cc}
m\\
k +  1\ \end{array})(a^{k +  1} - 1)a^{m - k - 1}b\}c^{m - k - 1}(a - 1)^{ - m +  1}\prod
\limits _{t=2}^{m - 1}(a^{t} - 1)\\
 &=& \frac {(m - 1)_a!} {c(1 - a)} ( (c - 1)^m (1 - ba^m c)  +
 (ac - 1)^m (cb  - 1)).
 \end {eqnarray*}

 {\rm (vii)} It follows from Part {(v)}.

 {\rm (viii)} It follows from Part {(v)} and Part {(vi)}.

{ \rm (ix)} We can show this by induction on $m$.

 {\rm (x)} It follows from {\rm (v)}.
 \hfill $\Box$

\vskip.1in
We use $[x_{i_{1}}\cdots x_{i_{m}}]$ instead of $[x_{i_{1}}[x_{i_{2}}[x_{i_{3}}\cdots x_{i_{m}}]]\cdots]$,   $1<m\leq n$. Similarly we define $[x_{i_{1}}\cdots x_{i_{m}}]^ -  := [x_{i_{1}}[x_{i_{2}}[x_{i_{3}}\cdots x_{i_{m}}]^ - ]^ - \cdots]^ - $.

\begin{Lemma}\label{3.2} If $<y_i,   v > =0$ and $v$ is a homogeneous element,   then

{\rm (i)} $<y_{i},  [x_i,   v]>=(p_{iv}^{ - 1} - p_{vi})v$,   $<y_{i}y_{j},  [x_i ,   v]>=(p_{iv}^{ - 1}p_{ij} - p_{vi}p_{ji}^{ - 1})<y_{j},  v>,  j\neq i$;

{\rm (ii)} $<y_{i},  [x_i,   v]^ - >=(1 - p_{iv}^{ - 1})v$,   $<y_{i}y_{j},  [x_i ,   v]^ - >=(p_{ji}^{ - 1} - p_{iv}^{ - 1}p_{ij})<y_{j},  v>,  j\neq i$.
\end {Lemma}

\begin{Lemma}\label{3.3}  Assume $ p_{i_s,   i_{t}} \begin {cases} \not=1,& t=s+1\\
=1 , & s +1 < t  \end  {cases}. $
 If $i_1,   i_2,   \cdots,   i_m$ is different each other and  $u_{sk}^ -  :=
[x_{i_s},   x_{i_{s +  1}} ,\cdots,   x_{i_k}]^ - $,    then

{\rm (i)} $u_{sk} ^ - \not=0$ for any $1\le s <k\le m$.

{\rm (ii)} $[u_{st}^ -  ,   u_{s +  1,   t} ^ - ]^ - \not=0$ for $1< s +  1 <t\le m$ and  $p _{i_k, i_k}\not= -1$ for $ s+ 1 \le k \le t$.
\end {Lemma}

\noindent{\bf Proof.} By Lemma \ref {3.2}(ii),   $<y_{i_s},   u_{sk}^ - > = (1 -  p_{i_s,    i_{s +  1}}^{ - 1}) u_{s +  1,   k} ^ -   $. We can show that $u_{sk} ^ - \not=0$ by induction on $\mid  k - s\mid$.
It is clear that
  $<y_{i_s},   [u_{st}^ - ,   u_{s +  1,   t}^ - ]^ - > = (1 -  p_{i_s,    i_{s +  1}}^{ - 1})^2 (u_{s +  1,   t} ^ - )^2$,   $<y_{i_s}^2,   (u_{st}^ - )^2> = (1+ p_{i_s, i_s} ^{-1}) p_{i_{s},   i_{ s +  1}}^{ - 1} (1 -  p_{i_{s},   i_{ s +  1}}^{ - 1})^2 u_{s +  1,   t}^2$ and $<y_{i_s}^2,   (u_{s,   s +  1}^ - )^2> = (1+ p_{i_s, i_s} ^{-1}) p_{i_{s},   i_{ s +  1}}^{ - 1} (1 -  p_{i_{s},   i_{s +  1}}^{ - 1})^2 x_{s +  1}^2.$  Consequently,   $<y_{i_{t - 1}}^2 \cdots y_{i_{s +  1}}^2y_{i_s},   [u_{st}^ - ,   u_{s +  1,   t}^ - ]^ - > = \{ (1+ p_{i_{s+1}, i_{s+1}} ^{-1})(1+ p_{i_{s+2}, i_{s+2}} ^{-1}) $ $\cdots  $  $ (1 + p_{i_{t-1},   i_{{t-1}}}^{ - 1} ) \} \{ p_{i_{s+ 1},   i_{s+2}}^{ - 1}  \cdots p_{i_{t - 1},   i_t}^{ - 1} \}\{ (1 -  p_{i_s,   i_{s +  1}}^{ - 1})^2(1 -  p_{i_{s +  1},   i_{s +  2}}^{ - 1})^2 \cdots  (1 -  p_{i_{t - 1},  i_t}^{ - 1}) ^2 \}x_{i_t}^2 $.

\hfill $\Box$

\begin{Lemma}\label{3.4} {\rm  (i)} Assume $ p_{i_s,   i_{t}} \begin {cases} \not=1,& t=s+1\\
=1 , & s +1 < t  \end  {cases}. $  If $i_1,   i_2,   \cdots i_m$ is different each other,   then
\begin {eqnarray} \label {e21.4.1}\dim \mathfrak L^ - (V) \ge  \sum \limits _{k=2}^{m}
(N_{i_k}  - 2)   +   C_m^2    +  n  + \mid C \mid, \end {eqnarray}   where $C:= \{ [u_{s,  t} ^ - ,   u _{{s +  1},   t}^ - ]^ -  \mid  p_{i_k, i_k} \not= -1,$ $s+1 \le k \le t  \}$ and   $ N_{i_k} := {\rm ord } (p_{i_k,   i_k})$ when $p_{i_k, i_k} \not=1$; $ N_{i_k} := \infty$  when $p_{i_k, i_k} =1$.

{\rm (ii)} Assume $p_{i_s,   i_{s +  1}} \not=1,  p_{i_{s+1},   i_{s }} \not=1, $
$p_{i_s,   i_{t}} =1,   $  for any $s  +  1 < t \le m$ and $1\le s < m$. If $i_1,   i_2,   \cdots i_m$ is different each other  then
\begin {eqnarray} \label {e21.4.2}\dim \mathfrak L^ - (V) \ge \sum \limits _{k=1}^{m - 1}
(N_{i_k}  - 2)  +  \sum \limits _{k=2}^{m}
(N_{i_k}  - 2)   +   C_m^2  +  n + \mid C \mid,\end {eqnarray}   where $C:= \{ [u_{s,  t} ^ - ,   u _{{s +  1},   t}^ - ]^ -  \mid  p_{i_k, i_k} \not= -1,$ $s+1 \le k \le t \}$ $\cup \{ [u_{s,  t} ^ - ,   u _{{s},   t +  1}^ - ]^ -  \mid  1\le  s<t < m, p_{i_k, i_k} \not= -1,$ $s \le k \le t \}$ and   $ N_{i_k} := {\rm ord } (p_{i_k,   i_k})$ when $p_{i_k, i_k} \not=1$; $ N_{i_k} := \infty$  when $p_{i_k, i_k} =1$.

\end {Lemma}
\noindent {\bf Proof.} {\rm (i)}
 Let   $A:=\{  \overline{l} _{i_{s+1}} ^k(i_s)^ -  \mid \ 2\le k < N_{i_{s+1}},      $ $ 1\le s <   m \}$,   $B:= \{  u_{s,   t} ^ - \mid  1\le s<t \le m\}$;
 $E:= \{ x_j \mid 1\le j \le n\}$.
By Lemma \ref {3.1} and Lemma \ref {3.3},   $A\cup B\cup C \cup E$ is linearly independent. It is clear that  $\mid A\mid =   \sum \limits _{k=2}^{m}
(N_{i_k}  - 2);$     $\mid B \mid =   C_m^2$.
Thus $\dim \mathfrak L^ - (V) \ge  \mid A\cup B\cup C \cup E \mid.$

{\rm (ii)} Let   $A:=\{  \overline{l} _{i_{t}} ^k(i_s)^ -  \mid \ 2\le k < N_{i_t},      $
  $  \mid t- s \mid =1$, $ 1\le s, t \le   m \}$,   $B:= \{  u_{s,   t} ^ - \mid  1\le s<t \le m\}$;
 $E:= \{ x_j \mid 1\le j \le n\}$.
By Lemma \ref {3.1} and Lemma \ref {3.3},   $A\cup B\cup C \cup E$ is linearly independent. It is clear that  $\mid A\mid = \sum \limits _{k=1}^{m - 1}
(N_{i_k}  - 2)  +  \sum \limits _{k=2}^{m}
(N_{i_k}  - 2);$     $\mid B \mid =   C_m^2$.
Thus $\dim \mathfrak L^ - (V) \ge  \mid A\cup B\cup C \cup E \mid.$
\hfill $\Box$

\begin{Theorem}\label{3.5} Assume that  $V$ is a connected  finite Cartan type  with   Cartan matrix $(a_{ij})_{n\times n}$. Let $N:={\rm ord } (q_{11})$ and $1\not= q \in F;$ $q_{ij} = q_{ji}$, $1\le i, j \le n.$

\vskip.2in
\noindent {\rm (i)}
  $A_n,$   $n>1$:    $\begin{picture}(100,    15) \put(27,   1){\makebox(0,
0)[t]{$\bullet$}} \put(60,    1){\makebox(0,   0)[t]{$\bullet$}}
\put(93,    1){\makebox(0,    0)[t]{$\bullet$}} \put(28,   -1){\line(1,
0){33}} \put(61,    -1){\line(1,    0){30}} \put(20,    - 15){1} \put(58,
- 15){2} \put(83,    - 15){3} \put(110,   0){$\dots$}\put(130,
1){\makebox(0,   0)[t]{$\bullet$}}\put(163,   1){\makebox(0,
0)[t]{$\bullet$}}\put(129,   -1){\line(1,   0){33}} \put(130,   - 15){n-1}
\put(163,  - 15){n}
\end{picture}$ \\

\vskip.1in
where $q_{ij} q_{ji} =q^{-1}$, $q_{ii} =q$ for $1\le i, j \le n$.

Case 1. $q^2 = 1.$  Then $\dim \mathfrak L^-(V)  \ge C_n^2 + n.$

Case 2. $q^2 \not= 1.$ Then
$\dim \mathfrak L^-(V) \ge 2 (n-1)
(N-2)  + 3C_n^2 - 2(n-1) + n  =2 (n-1)
(N-2)  + 3 C_n^2 -n +2$.

\vskip.2in
\noindent
{\rm (ii)}  $B_n$, $n>2$:

\ \ \ \ \ \ \ \ \  \ \ \
$
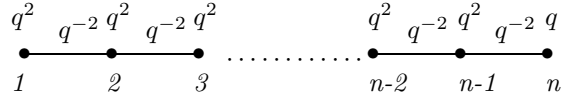
\begin{picture}(100,     15)
\put(27,     1){\makebox(0,    0)[t]{$\bullet$}}
\put(60,     1){\makebox(0,     0)[t]{$\bullet$}}
\put(93,    1){\makebox(0,    0)[t]{$\bullet$}}
\put(159,     1){\makebox(0,     0)[t]{$\bullet$}}
\put(192,    1){\makebox(0,     0)[t]{$\bullet$}}
\put(225,    1){\makebox(0,    0)[t]{$\bullet$}}
\put(28,     -1){\line(1,     0){33}}
\put(61,     -1){\line(1,     0){30}}
\put(130,    -1){\makebox(0,    0)[t]{$\cdots\cdots\cdots\cdots$}}
\put(160,    -1){\line(1,     0){30}}
\put(193,     -1){\line(1,     0){30}}
\put(22,    -15){1}
\put(58,     -15){2}
\put(91,     -15){3}
\put(157,     -15){n-2}
\put(191,     -15){n-1}
\put(224,     -15){n}
\put(22,    10){$q^{2}$}
\put(58,     10){$q^{2}$}
\put(91,     10){$q^{2}$}
\put(157,     10){$q^{2}$}
\put(191,     10){$q^{2}$}
\put(224,     10){$q$}
\put(40,     5){$q^{-2}$}
\put(73,     5){$q^{-2}$}
\put(172,    5){$q^{-2}$}
\put(205,     5){$q^{-2}$}
\end{picture}$\\

\vskip.1in
$N_i = N = {\rm ord }(q^2),$ $N_n = {\rm ord }(q),$ $1\le i \le n-1.$

 Case 1. $q^4 = 1$ and $q^2 \not=1.$  Then $\dim \mathfrak L^-(V)  \ge (N_{n}  - 2)   + C_n^2  + n = N_n -2 + C_n^2 + n$.

Case 2. $q^4 \not= 1.$  Then $\dim \mathfrak L^-(V) \ge  \sum \limits _{k=1}^{n - 1}
(N_{k}  - 2)  +  \sum \limits _{k=2}^{n}
(N_{k}  - 2)   + 3 C_n^2 - 2(n-1) + n = 2(n-1) (N-2) + (N_n -2) - (N-2) + 3 C _n^2-n +2.$

\vskip.2in
\noindent
{\rm (iii)} $C_n$, $n>2$:

\ \ \ \ \ \ \ \ \  \ \ \
$
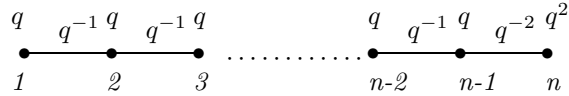
\begin{picture}(100,     15)
\put(27,     1){\makebox(0,    0)[t]{$\bullet$}}
\put(60,     1){\makebox(0,     0)[t]{$\bullet$}}
\put(93,    1){\makebox(0,    0)[t]{$\bullet$}}
\put(159,     1){\makebox(0,     0)[t]{$\bullet$}}
\put(192,    1){\makebox(0,     0)[t]{$\bullet$}}
\put(225,    1){\makebox(0,    0)[t]{$\bullet$}}
\put(28,     -1){\line(1,     0){33}}
\put(61,     -1){\line(1,     0){30}}
\put(130,    -1){\makebox(0,    0)[t]{$\cdots\cdots\cdots\cdots$}}
\put(160,    -1){\line(1,     0){30}}
\put(193,     -1){\line(1,     0){30}}
\put(22,    -15){1}
\put(58,     -15){2}
\put(91,     -15){3}
\put(157,     -15){n-2}
\put(191,     -15){n-1}
\put(224,     -15){n}
\put(22,    10){$q$}
\put(58,     10){$q$}
\put(91,     10){$q$}
\put(157,     10){$q$}
\put(191,     10){$q$}
\put(224,     10){$q^2$}
\put(40,     5){$q^{-1}$}
\put(73,     5){$q^{-1}$}
\put(172,    5){$q^{-1}$}
\put(205,     5){$q^{-2}$}
\end{picture}$\\

\vskip.1in
$N_i = N = {\rm ord }(q),$ $N_n = {\rm ord }(q^2),$ $1\le i \le n-1.$

   Case 1. $q^4 = 1$ and $q^2 \not=1$.  Then $\dim \mathfrak L^ - (V) \ge    \sum \limits _{k=1}^{n - 2}
(N_{k}  - 2)  +  \sum \limits _{k=2}^{n-1}
(N_{k}  - 2)   + 3 C_{n-1}^2 - 2(n-2) + n-1 +\{  N_{n-1}-2 +(n-1) + (n-2) +1 \} =
2(n-2) (N-2) + (N-2) + 3C_{n-1}^2 + n-1$.

Case 2. $q^4 \not= 1.$
$\dim \mathfrak L^-(V) \ge  \sum \limits _{k=1}^{n - 1}
(N_{k}  - 2)  +  \sum \limits _{k=2}^{n}
(N_{k}  - 2)  + 3 C_n^2 - 2(n-1) + n   = 2(n-1) (N-2) + (N_n -2) - (N-2) + 3 C _n^2 -n+2.$

\vskip.2in
{\rm (iv)} $D_n (n>3)$: generalized Dynkin diagram:\\  \\
$
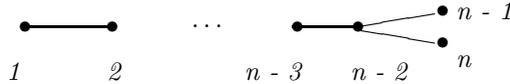
\begin{picture}(200,    15)  \put(130,    7){\makebox(0,
0)[t]{$\bullet$}} \put(153,    7){\makebox(0,    0)[t]{$\bullet$}}
 \put(185,    1){\makebox(0,    0)[t]{$\bullet$}}
\put(185,   13){\makebox(0,    0)[t]{$\bullet$}}  \put(110,    - 15){n - 3}
\put(150,    - 15){n - 2} \put(190,     - 10){n} \put(190,    8){n - 1}
\put(152,   5){\line(6,    1){30}} \put(182,     - 1){\line( - 6,    1){30}}

 \put(128,    5){\line(1,    0){25}}

 \put(27,   7){\makebox(0,   0)[t]{$\bullet$}} \put(60,    7){\makebox(0,
0)[t]{$\bullet$}} \put(28,   5){\line(1,   0){33}}  \put(20,   - 15){1}
\put(58,   - 15){2}  \put(90,   5){$\dots$}

\end{picture}$\\

\vskip.1in
where $q_{ij} q_{ji} =q^{-1}$, $q_{ii} =q$ for $1\le i, j \le n$.

Case 1.  $q^2 = 1.$  Then $\dim \mathfrak L^-(V)  \ge  C_{n-1}^ 2 +2n -1.$

Case 2.  $q ^2\not= 1 .$
Then $\dim \mathfrak L^-(V) \ge (2 (n-1-1)
(N-2)  + 3 C_{n-1}^2 - 2(n-1-1) + n-1) + ( 2 (N-2) +  (n-1 ) + (n-3) + (n-3) +2 +1) = 2 (n-2) (N-2) +3 C_{n-1}^2 + 2n +2N-5 $.

\vskip.2in
{\rm (v)}  $E_8$:
  \\
  \\
  $
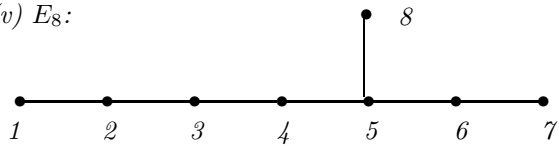
\begin{picture}(100,    15) \put(158,   34){\makebox(0,
0)[t]{$\bullet$}} \put(27,    1){\makebox(0,    0)[t]{$\bullet$}}
\put(60,    1){\makebox(0,    0)[t]{$\bullet$}} \put(93,
1){\makebox(0,   0)[t]{$\bullet$}} \put(126,    1){\makebox(0,
0)[t]{$\bullet$}} \put(159,    1){\makebox(0,    0)[t]{$\bullet$}}
\put(192,   1){\makebox(0,    0)[t]{$\bullet$}} \put(225,
1){\makebox(0,   0)[t]{$\bullet$}} \put(28,     - 1){\line(1,    0){33}}
\put(157,   1){\line(0,    1){30}} \put(61,     - 1){\line(1,    0){30}}
\put(94,    - 1){\line(1,    0){30}} \put(127,     - 1){\line(1,    0){30}}
\put(160,    - 1){\line(1,    0){30}} \put(193,     - 1){\line(1,    0){30}}
\put(22,    - 15){1} \put(58,     - 15){2} \put(91,     - 15){3} \put(124,
 - 15){4} \put(157,     - 15){5} \put(191,     - 15){6} \put(224,     - 15){7}
\put(170,   28){8}
\end{picture}$\\
\vskip.1in
where $q_{ij} q_{ji} =q^{-1}$, $q_{ii} =q$ for $1\le i, j \le n$.

Case 1.  $q^2 = 1.$  Then $\dim \mathfrak L^-(V)  \ge C_{8-1}^2 + 8 +  7 = C_7^2 +15.$

Case 2.  $q ^2\not= 1 .$
Then
$\dim \mathfrak L^-(V)
\ge ( 2 (7-1)
(N-2)   + 3 C_7^2 - 2(7-1) +7 ) + (2(N-2)+7+4 + 4+4+1 ) =   14N+50. $

\vskip.2in
\noindent
{\rm (vi)} $E_7$.
Case 1.  $q^2 = 1.$  Then $\dim \mathfrak L^-(V)  \ge C_{7-1}^2 + 7 +  6 = C_6^2 + 13.$

Case 2.  $q ^2\not= 1 .$
Then $\dim \mathfrak L^-(V)
\ge ( 2 (6-1)
(N-2)   + 3C_6^2 - 2(6-1) +6 ) + (2(N-2)+6+3+3+4+1 ) = 12 N +34. $

\vskip.2in
\noindent
{\rm (vii)} $E_6.$
Case 1.  $q^2 = 1.$  Then $\dim \mathfrak L^-(V)  \ge C_{6-1}^2 + 6 +  5.$

Case 2.  $q ^2\not= 1 .$
 Then $\dim \mathfrak L^-(V)
\ge ( 2 (5-1)
(N-2)   + 3 C_5^2 - 2(5-1) +5 ) + (2(N-2)+5+2+2+4+1 ) = 10N +21.$

\vskip.2in
\noindent
{\rm (viii)} $F_4$.
$\begin{picture}(100,     15)
\put(27,     1){\makebox(0,    0)[t]{$\bullet$}}
\put(60,     1){\makebox(0,     0)[t]{$\bullet$}}
\put(93,    1){\makebox(0,    0)[t]{$\bullet$}}
\put(126,     1){\makebox(0,   0)[t]{$\bullet$}}
\put(28,     -1){\line(1,     0){33}}
\put(61,     -1){\line(1,     0){30}}
\put(94,    -1){\line(1,     0){30}}
\put(22,    -15){1}
\put(58,     -15){2}
\put(91,     -15){3}
\put(124,    -15){4}
\put(22,    10){$q^2$}
\put(58,     10){$q^2$}
\put(91,     10){$q$}
\put(124,     10){$q$}
\put(40,     5){$q^{-2}$}
\put(73,     5){$q^{-2}$}
\put(106,    5){$q^{-1}$}
\end{picture}$\\

\vskip.1in
Case 1. $q^4 = 1$ and $q^2 \not= 1.$  Then $\dim \mathfrak L^-(V)  \ge  2(N_3 -2) + (N_4-2)   + C_4^2  + 4 + 4$.

Case 2. $q^4 \not= 1.$
Then $\dim \mathfrak L^-(V) \ge \sum \limits _{k=1}^{4 - 1}
(N_{k}  - 2)  +  \sum \limits _{k=2}^{4}
(N_{k}  - 2)  + 3 C_4^2 - 2(4-1) + 4 .$

\vskip.2in
\noindent
{\rm (ix)} $G_2$:
\ \ \ \ \ \ \ \ \  \ \ \
$\begin{picture}(100,     15)
\put(27,     1){\makebox(0,    0)[t]{$\bullet$}}
\put(60,     1){\makebox(0,     0)[t]{$\bullet$}}
\put(28,     -1){\line(1,     0){33}}
\put(22,    -15){1}
\put(58,     -15){2}
\put(22,    10){$q$}
\put(58,     10){$q^3$}
\put(40,     5){$q^{-3}$}
\end{picture}$\\
\vskip.1in  Case 1. $q^2 = 1.$  Then $\dim \mathfrak L^-(V)  \ge  1+2. $

Case 2. $q^2 \not= 1.$
 Then $\dim \mathfrak L^-(V) \ge
 \sum \limits _{k=1}^{2 - 1}
(N_{k}  - 2)  +  \sum \limits _{k=2}^{2}
(N_{k}  - 2)     + 3  .$

\end {Theorem}
\noindent {\bf Proof.} { \rm (i)} $N_i = N = {\rm ord }(q).$  Case 1.
Let $A: = \emptyset;$ $B:= \{  u_{s,   t} ^ - \mid  1\le s<t \le n\}$; $C: = \emptyset.$

Case 2. Let   $A:=\{  \overline{l} _{{t}} ^k(s)^ -  \mid \ 2\le k < N_{t},      $
  $  \mid t- s \mid =1$, $ 1\le s, t \le   n \}$,   $B:= \{  u_{s,   t} ^ - \mid  1\le s<t \le n\}$;
$C:= \{ [u_{s,  t} ^ - ,   u _{{s +  1},   t}^ - ]^ -  \mid  $ $1< s +  1<t \le n\}$ $\cup \{ [u_{s,  t} ^ - ,   u _{{s},   t +  1}^ - ]^ -  \mid  1\le  s<t < n \}$; $E:= \{ x_j \mid 1\le j \le n\}$.
By Lemma \ref {3.1} and Lemma \ref {3.3},   $A\cup B\cup C \cup E$ is linearly independent.
Thus $\dim \mathfrak L^ - (V) \ge  \mid A\cup B\cup C \cup E \mid.$

\vskip.2in

{\rm (ii)} $N_i = N = {\rm ord }(q^2),$ $N_n = {\rm ord }(q),$ $1\le i \le n-1.$

Case 1. Let $A: = \{  \overline{l} _{{n}} ^k(n-1)^ -  \mid \ 2\le k < N_{n}\},      $
   ; $B:= \{  u_{s,   t} ^ - \mid  1\le s<t \le n\}$;
    $C:= \emptyset$; $E:= \{ x_j \mid 1\le j \le n\}$. Consequently,  $\dim \mathfrak L^-(V)  \ge  (N_n -2)   + C_n^2   + n$.

   Case 2. It  follows from Lemma \ref {3.4} (ii) and the proof of  Part (i).

\vskip.2in

   (iii) $N_i = N = {\rm ord }(q),$ $N_n = {\rm ord }(q^2),$ $1\le i \le n-1.$

   Case 1. Let   $A:=\{  \overline{l} _{{t}} ^k(s)^ -  \mid \ 2\le k < N_{t},      $
  $  \mid t- s \mid =1$, $ 1\le s, t \le   n-1 \}\cup \{ \overline{l} _{{n-1}} ^k(n)^ -  \mid \ 2\le k < N_{n-1}  \}$,   $B:= \{  u_{s,   t} ^ - \mid  1\le s<t \le n\}$;
$C:= \{ [u_{s,  t} ^ - ,   u _{{s +  1},   t}^ - ]^ -  \mid  $ $1< s +  1<t < n\}$ $\cup \{ [u_{s,  t} ^ - ,   u _{{s},   t +  1}^ - ]^ -  \mid  1\le  s<t < n, \}$; $E:= \{ x_j \mid 1\le j \le n\}$.
By Lemma \ref {3.1} and Lemma \ref {3.3},   $A\cup B\cup C \cup E$ is linearly independent.
Thus $\dim \mathfrak L^ - (V) \ge  \mid A\cup B\cup C \cup E \mid=  \sum \limits _{k=1}^{n - 2}
(N_{k}  - 2)  +  \sum \limits _{k=2}^{n-1}
(N_{k}  - 2)   + 2 C_{n-1}^2 - (n-3) + n-1 +\{  N_{n-1}-2 +(n-1) + (n-2) +1 \}$

   Case 2. It  follows from Lemma \ref {3.4} (ii) and the proof of  Part (i).

\vskip.2in

{\rm (iv)} Let $A_1:=\{  \overline{l} _{s} ^k(t)^ -  \mid \ 2\le k < N,    \mid t -  s \mid =1,   $ $ 1\le s; t \le n - 1\}$,   $B_1:= \{  u_{s,   t} ^ - \mid  1\le s<t \le n - 1\}$;
$C_1:= \{ [u_{s,  t} ^ - ,   u _{{s +  1},   t}^ - ]^ -  \mid  1< s +  1<t \le n - 1\}\cup \{ [u_{s,  t} ^ - ,   u _{{s},   t +  1}^ - ]^ -  \mid  1\le  s<t < n - 1\}$; $E_1:= \{ x_j \mid 1\le j \le n - 1\}$;
$A_2:=\{  \overline{l} _{s} ^k(t)^ -  \mid  s =n - 2, t=n,$ or $ s= n, t= n - 2\}$,   $B_2:= \{  [x_s, x_{s +  1}, \cdots, x_{n-2}, x_n ] ^ - , [x_{n - 1}, x_ {n - 2}, x_n]^ -   \mid  1\le s \le n - 2\}$;
$C_2:= \{ [[x_s, x_{s +  1}, \cdots, x_{n - 2}, x_n ] ^ - ,  [x_{s +  1}, x_{s +  2}, \cdots, x_{n - 2}, x_n ] ^ - ]^ - \mid  1\le s<n - 2 \} $ $\cup $ $\{ [[x_s, x_{s +  1}, \cdots, x_{n - 2}, x_n ] ^ - , $ $[x_{s}, x_{s +  1}, \cdots, x_{n - 2} ] ^ - ]^ -  $ $\mid  1\le s < n - 2 \} $ $\cup \{ [[ x_{n - 1}, x_{n - 2}, x_n]^ -  , [  x_{n -2},  x_n]^ - ]^ - \} $ $\cup \{ [[ x_{n - 1}, x_{n - 2}, x_n]^ -  , [  x_{n - 1},  x_{n - 2}]^ - ]^ - \}; $ $E_2:= \{ x_n \}$.

\vskip.2in

 { \rm (v)}
Let $A_1:=\{  \overline{l} _{s} ^k(t)^ -  \mid \ 2\le k < N,    \mid t -  s \mid =1,   $ $ 1\le s; t \le 7\}$,   $B_1:= \{  u_{s,   t} ^ - \mid  1\le s<t \le 7\}$;
$C_1:= \{ [u_{s,  t} ^ - ,   u _{{s +  1},   t}^ - ]^ -  \mid  1< s +  1<t \le 7\}\cup \{ [u_{s,  t} ^ - ,   u _{{s},   t +  1}^ - ]^ -  \mid  1\le  s<t < 7\}$; $E_1:= \{ x_j \mid 1\le j \le 7\}$;
$A_2:=\{  \overline{l} _{s} ^k(t)^ -  \mid  s =5, t=8,$ or $ s= 8, t= 5\}$,   $B_2:= \{  [x_s, x_{s +  1}, \cdots, x_5, x_8 ] ^ - , [x_7, x_6, x_5, x_8]^ - , [ x_6, x_5, x_8]^ -  \mid  1\le s \le 5\}$; \\
$C_2:= \{ [[x_s, x_{s +  1}, \cdots, x_5, x_8 ] ^ - ,  [x_{s +  1}, x_{s +  2}, \cdots, x_5, x_8 ] ^ - ]^ - \mid  1\le s<5 \} $ $\cup $ $\{ [[x_s, x_{s +  1}, \cdots, x_5, x_8 ] ^ - , $ $[x_{s}, x_{s +  1}, \cdots, x_5 ] ^ - ]^ -  $ $\mid  1\le s< 5 \} $ $\cup \{ [[ x_7, x_6, x_5,  x_8]^ -  , [  x_6, x_5,  x_8]^ - ]^ - ; $  $ [[ x_6,  x_5,  x_8]^ -  , [   x_5,  x_8]^ - ]^ - ; $ $[[ x_7, x_6, x_5,  x_8]^ -  , [ x_7,  x_6, x_5]^ - ]^ - ; $ $\ \ [[ x_6,  x_5,  x_8]^ -  , [ x_6,   x_5]^ - ]^ -  \}$; $E_2:= \{ x_j \mid 8\le j \le 8\}$.

\vskip.2in

 {\rm (vi)} and {\rm (vii)} They are similar to the proof of {\rm (v)}.

\vskip.2in

  {\rm (viii)} Case 1. Let $A: = \{  \overline{l} _{{4}} ^k(3)^ - , \overline{l} _{{3}} ^j(4)^ - , \overline{l} _{{3}} ^j(2)^ -\mid \ 2\le k < N_{4}, 2\le j < N_{3}\};     $
   ; $B:= \{  u_{s,   t} ^ - \mid  1\le s<t \le 4\}$;
    $C:= \{ [u_{1,  4} ^ - ,   u _{{1 +  1},   4}^ - ]^ -;  [u_{2,  4} ^ - ,   u _{{2 +  1},   4}^ - ]^ - ; [u_{1,  4} ^ - ,    u _{{  1},   3}^ - ]^ - ; [u_{2,  4} ^ - ,   u _{{2},   3}^ - ]^ -   \}$; $E:= \{ x_j \mid 1\le j \le 4\}$. Consequently, $\dim \mathfrak L^-(V)  \ge  2(N_2 -2) + (N_4-2)   + C_4^2  + 4 + 4$.

   Case 2. It  follows from Lemma \ref {3.4} (ii) and the proof of  Part (i).

\vskip.2in

{\rm (ix)} Case 1. Let $A: = \emptyset; $  $B:= \{  u_{1,   2} ^ -\}$;
    $E:= \{ x_j \mid 1\le j \le 2\}$. Consequently, $\dim \mathfrak L^-(V)  \ge  1+2. $

   Case 2. It  follows from Lemma \ref {3.4} (ii) and the proof of  Part (i).
 \hfill $\Box$
\vskip.1in
We remark that it is possible that the dimensions of Nichols Lie algebras of two twisting equivalent braided vector spaces are different.

\begin {Example} \label {3.6} Let  $V, V'$ and $ V''$ are three braided vector spaces with braiding matrixes  $(q_{ij})_{2\times 2}$, $(q_{ij}')_{2\times 2}$ and $(q_{ij}'')_{2\times 2}$, respectively. Assume $q_{1,  2} q_{2,  1} =1$,   $q_{1,  1} = q_{2,  2}= - 1$, $q_{ii} = q_{ii}'= q_{ii}''$, $q_{ij} q_{ji} = q_{ij}' q_{ji}'= q_{ij} ''q_{ji}''$ for $i, j=1,2$. Then

{\rm (i)} $\dim \mathfrak B(V) = 4.$

{\rm (ii)} $ [x_1,   x_2]^ -  =0$ and $\dim \mathfrak L^ - (V') = 2$ when  $q_{1,  2}' =1$.

{\rm (iii)}   $ [x_1,   x_2]^ -  \not= 0$ and $\dim \mathfrak L^ - (V'') = 3$ when  $q_{1,  2} '' \not=1$.

{\rm (iv)}  $V'$ and $V''$ are twisting equivalent. \end {Example}

\begin {Lemma} \label {3.7} {\rm (i)} If there exists $1\le i \le n$ such that ${\rm ord }(p_{ii}) =1 $ or $\infty,$ then there exists  $V'$, which is twisting equivalent to $V$, such that  $ \dim (\mathfrak L^-(V')) = \infty.$

{\rm (ii)} If $p_{11}=p_{22} =-1$, then
 $< (y_2y_1)^ k , x_1 (x_2x_1) ^m> =
 \prod\limits _{j=0} ^{k-1}(1 - (p_{12}p_{21})^{-m-j}) x_1 (x_2x_1)^{m-k}$  for $1 \le k \le m.$

{\rm (iii)} If $p _{11}  = p_{22} =-1$ and ${\rm ord } (p_{12} p _{21}) = \infty,$  then
$ \bar l _ {[x_2, x_1]^- } ^m(x_1)^- \not= 0 $ for any $m \in \mathbb N,$ which implies
$\dim (\mathfrak L^-(V)) = \infty.$

{\rm (iv)} If  the generalized Dynkin diagram of  $V$ is a simple chain with length $d$ and there exists $1\le i < d$ such that ${ \rm ord } (p_{i, i+1} p_{i+1, i}) = \infty$, then  there exists $1\le j \le n$ such that ${\rm ord } (p_{jj}) = \infty$ or $\dim ( \mathfrak L^-(V)) = \infty.$

{\rm (v)} If $\Delta (\mathfrak B(V))$ is an arithmetic root system and  there exists $u\in D(V)$ such that ${\rm ord } (p_{uu}) = \infty$, then there exists $1\le i \le n$ such that ${\rm ord } (p_{ii}) = \infty$ or $\dim ( \mathfrak L^-(V)) = \infty.$

\end {Lemma}
\noindent{\bf Proof.}  { \rm (i)} It follows from Lemma \ref {3.1} { \rm (x)}.

{\rm (ii)} It can be proved by induction on $m$ that $<y_1,  x_1 (x_2x_1) ^m  > =  (x_2x_1) ^m +  p _{1, 12}^{-m} (x_1x_2)^m$ and $<y_2y_1,  x_1 (x_2x_1) ^m  > =   (1 - ( p_{12} p_{21})^ {-m}) x_1 (x_2x_1) ^{m-1}. $ It can be proved by induction on $k$ that
$< (y_2y_1)^ k , x_1 (x_2x_1) ^m > =
 \prod \limits_{j=0} ^{k-1}(1 - (p_{12}p_{21})^{-m-j}) x_1 (x_2x_1)^{m-k}.$

{\rm  (iii)}  It follows from {\rm (ii)}.

{\rm (iv) }  If $1\le i <d$ and $q_{ii} \not= -1$, then   ${\rm ord }(q_{ii}) = \infty$ by Definition 1 in 
\cite {He06a}. Similarly, $q_{i+1, i+1} \not= -1$, then   ${\rm ord }(q_{i, i}) = \infty$.
If   $q_{ii} = -1$ and $q_{i+1, i+1} = -1$,  then we complete the proof by Part {\rm (ii)}.

{\rm (v) }   If there exists $u \in D(V)$  such that ${\rm ord }  (p_{uu}) =\infty$, then there exists a $1\le i \le n$  such that ${\rm ord }  (p_{ii}) = \infty $  except the cases of  Row 3, Diagram 2, Table 1;  Row 8, Diagram 2,  Table 2; Row 9, Diagram 4,  Table 2;  Row 10, Diagram 3,  Table 2, in \cite{He05}, and the cases of  Row 10, Diagram 6,  Appendix B;  Row 12, Diagram 5,  Appendix B. Row 1-10,   Appendix C, in \cite{He06a}.  By \rm {(iii)},
$\dim (\mathfrak L^-(V)) = \infty.$\hfill $\Box$

\begin {Theorem} \label {3.8} Assume that  $\Delta (\mathfrak B(V))$ is an arithmetic root system
and there does not exist any $m$-infinity element with $p_{u,u} \not= 1$ for any $u \in D(V)$. Then

{\rm (i)}   $\dim (\mathfrak B(V) ) = \infty$  if and only if there exists  $V'$,  which is twisting equivalent to $V$, such that $ \dim (\mathfrak L^ - (V')) = \infty.$

{\rm (ii)}   $\mathfrak B(V)$ is
finite-dimensional if and only if $\mathfrak L^-(V)$ is finite-dimensional when the input of  every vertex in the directed Dynkin diagram of $V$ is more than $0$, i.e. for any $1\le i \le n$, there exists $j$ such that  $p_{ji} \not= 1$.

\end {Theorem}

\noindent{\bf Proof.}   By \cite {He05} $\mathfrak B(V)$ is
finite-dimensional if and only if $1< {\rm ord } (p_{uu}) < \infty$ for any $u\in D(V)$.

{\rm (i)} The sufficiency is clear. The necessity.  By \cite {He05}  there exists a $u\in D(V)$ such that ${\rm ord }  (p_{uu}) =1$ or $\infty$. We complete the proof by Lemma \ref {3.7} {\rm (i), (iv), (v)}.

{\rm (ii)}  It follows from  Lemma \ref {3.1}(x) and the proof of Part {\rm (i)}. \hfill $\Box$

\begin {Corollary} \label {3.9} Assume that  $\Delta (\mathfrak B(V))$ is an arithmetic root system. Then
 $\dim (\mathfrak B(V) ) = \infty$  if and only if there exists  $V'$,  which is twisting equivalent to $V$, such that $ \dim (\mathfrak L^ - (V')) = \infty$ in the following three cases:

{\rm (i)}  $\dim V=2.$

{\rm (ii)}   $V$ is of finite Cartan type.

{\rm (iii)}  $V$ is a  Yetter-Drinfeld module over finite cyclic groups.

\end {Corollary}
\noindent{\bf Proof.} It is enough to show that  there does not exist any $m$-infinity element and  $p_{u,u} \not= 1$ for any $u \in D(V)$ in this three cases.

There are no $m$-infinity elements in this three cases by Theorem \ref {2.10}
and
 Lemma \ref {8'}. The second claim is obtained by
Appendix, Lemma 6.4 in \cite {WZZ15} and Theorem 3.3 in \cite {WZZ14}. \hfill $\Box$

\begin {Theorem} \label {3.10} If $V$ is a quantum linear space and $0\not=\mid p _{s, t} \mid  <1$  when  $t=s +  1$; $p_{st} = 1$  when $t\not= s +  1$ and $s< t, $ then
$\dim \mathfrak L^ -  (V) = N_1 \cdots N_n  -  (N_1  +   \cdots  +   N_n  - n+1) $ and
$ \mathfrak B (V)  =  \mathfrak L^ -  (V) \oplus span \{ x_i^ {m_i}  \mid  1< m_i \in \mathbb N \} \oplus F.$
\end {Theorem}

\noindent {\bf Proof.} We show this by following several steps.

 Let $W = x_1 ^{a_1} \cdots x_n ^{a_n}$ and $N(W) : = a_1  +   \cdots a_n$. Assume $a_{i_0}$ is the first non-zero element and $a_{j_0}$ is the lost non-zero element.

{\rm (i)} If $u =x_1 ^{a_1} \cdots x_r ^{a_r} $  and $v=x_{r +  1} ^{a_{r +  1}} \cdots x_n ^{a_n}$, then $uv = p_{u, v}vu \in \mathfrak L^ -  (V)$ and $p_{u, v} = \prod \limits_{i = 1, \cdots, r ;  j=r +  1, \cdots, n} p_{i, j}^{a_i a_j}$ when $u, v \in \mathfrak L^-(V)$.

{\rm (ii)} We use induction on $N(W)$ to show that $W \in \mathfrak L^ - (V)$ when $a_i \le 1$ for $1\le i \le n.$ It is clear for $N(W) =1.$  For $N(W)>1,$ Let $1\le i_0\le n$ such that $a_{i_0} $ is the first non-zero element. Let $u= x_{i_0}$, $v= x_{i_0 +  1}^{a_{i_0 +  1}} \cdots x_n^{a_n}.$ By induction we have $v\in \mathfrak L^ - (V)$.
Consequently, $W= uv \in \mathfrak L^ - (V)$ by {\rm (i)}.

{\rm (iii)} We use induction on $N(W)$ to show that $W \in \mathfrak L^ - (V)$ when there exist two non-zero elements in $ \{ a_1, a_2, \cdots, a_n\}$. When $N(W) =2$, this is Case {\rm (i)}.  For $N(W) =3,$ considering {\rm (i)} we can assume that there exist  only two  $a_{i_0} \not=0$  and $a_{j_0} \not=0$. If $a_{i_0 }>1$, then $a_{j_0}>1$.  Let  $u= x_{i_0}$, $v= x_{i_0} ^{a_{i_0}  - 1} x_{i_0 +  1}^{a_{i_0 +  1}} \cdots x_n^{a_n}.$  By induction, $ v \in  \mathfrak L^ - (V)  $. Consequently, $W \in \mathfrak L^ - (V)$. If $a_{i_0 }=1$, let  $u= x_{i_0}^{a_{i_0}} \cdots x_{j_0} ^{a_{j_0} - 1}$, $v= x_{j_0}.$  By induction, $ u\in  \mathfrak L^ - (V)  $. Consequently, $W= uv \in \mathfrak L^ - (V)$.
For $N(W)> 3, $ if $a_{i_0}=1$ and $a_{j_0} =1$, let $u= x_{i_0}$, $v= x_{i_0 +  1} ^{a_{i_0 +  1}} \cdots x_n^{a_n}$.  By induction, $ v\in  \mathfrak L^ - (V)  $. Consequently, $W= uv \in \mathfrak L^ - (V)$. If $a_{i_0 }>1$, let  $u= x_{i_0}$, $v= x_{i_0} ^{a_{i_0}  - 1} x_{i_0 +  1}^{a_{i_0 +  1}} \cdots x_n^{a_n}.$  By induction, $ v \in  \mathfrak L^ - (V)  $. Consequently, $W \in \mathfrak L^ - (V)$. If $a_{i_0 }=1$ and $a_{j_0} >1$,  let  $u= x_{i_0}^{a_{i_0}} \cdots x_{j_0} ^{a_{j_0} - 1}$, $v= x_{j_0}.$  By induction, $ u\in  \mathfrak L^ - (V)  $. Consequently, $W= uv \in \mathfrak L^ - (V)$.

{\rm (iv)}  By Theorem 2 in \cite  {Kh99b} or  Theorem 10 in \cite {He07}, $\{ x_1 ^{a_1} \cdots x_n ^{a_n} \mid   0 \le a_i < N_i, 1\le i \le n  \}$ is a basis of $\mathfrak B(V)$. Considering {\rm (i) - (iii)} and Lemma \ref {3.1} (ix) we complete the proof.

\hfill $\Box$

Let us remark that it is possible that $\dim \mathfrak B(V) = \infty$ and $\dim \mathfrak L^ - (V)< \infty.$
\begin {Example} \label {3.11} Assume $q_{1  2} = q_{2  1} =1$ and    $q_{1  1} = q_{2  2}= 1 $. Then
 $\dim \mathfrak B(V) = \infty$ and $\dim \mathfrak L^ - (V) = 2$.
 \end {Example}

\section {Two examples of Lie algebras without maximal solvable ideals}\label{solvable ideals}

In this section we find two examples of Lie algebras which have no maximal solvable ideals.

Let $R_k := \{ x \mid x  \hbox { is a primitive } k \hbox {-th unit root} \}$. Assume $(L, [\ ]^-)$ is a Lie algebra.
Let  $ D^{(0)}
(L) := L,   D^{(1)}(L) = [L,   L]^-,$ $D^{(k+1)}(L) = [D^{(k)}(L),   D^{(k)}(L)]^-$.
If there exists a natural number $m$ such that $D^{(m)}(L) =0,$ then  $L$ is called a solvable Lie algebra.

 For any $z_1, z_2, \cdots, z_{k}\in L$, let $\sigma  (z_1, z_2) := [z_1, z_2]^-$, $\sigma (z_1, z_2, z_3, z_4) := [\sigma (z_1, z_2), \sigma (z_3, z_4)]^- .$  $\cdots, $
$\sigma (z_1, z_2,\cdots, z_{2^k})$ $:= [\sigma ( z_1, z_2, \cdots, z_{2^{k-1}}), \sigma (z_{2^{k-1}+1},  \cdots,  z_{2^k})]^-.$

\begin {Lemma} \label{4.1}  Assume that $V_m$ is a  quantum vector space with $q_{i,i+1} \not=1$  and $q_{ij}=1$, $i+1 <j\le m, 1 \le i <m$. If $x_1, x_2, \cdots, x_{m} $
be the canonical basis of $V_{m}$ and  $k-s$ is a power of $2,$ then $<y_s, \sigma (x_s, \cdots x_{k-1})> \not=0$  for any $1 \le s < k \le   m+1.$

\end{Lemma}
\noindent {\bf Proof.}  Now we show that $< y_s, \sigma  (x_s, x_{s+1}, \cdots, x_{k-1})>$ $\not= 0$ by induction on $k-s$, where $k-s$ is a power of $2,$ for any $1 \le s < k \le   m+1.$  If $k-s =2,$ then $< y_s, \sigma  (x_s, x_{s+1}, \cdots, x_{k-1})> = (1- p_{s, s+1}^{-1}) x_{s+1}$ $\not= 0$. Assume $k-s >2.$ Then
\begin {eqnarray*} && <y_s,  \sigma  (x_s, x_{s+1}, \cdots, x_{k-1})> \\
&=&  (1-p^{-1}_{\frac {k-s} {2}, \frac {k-s} {2} +1})
<y_s, \sigma  (x_s, x_{s+1}, \cdots, x_{\frac {k-s}{2}  } )
> \sigma  (x_{\frac {k-s}{2} +1}, x_{\frac {k-s}{2}+2}, \cdots,x_{k-1} ).
  \end {eqnarray*}
  By induction assumption, $ <y_s, $ $\sigma  (x_s, x_{s+1}, $ $\cdots, $ $x_{\frac {k-s}{2}  } )
>\not=0$ and $ \sigma  (x_{\frac {k-s}{2} +1}, x_{\frac {k-s}{2}+2}, \cdots, x_{k-1} ) \not=0. $ Considering the restricted  PBW basis of Nichols algebras (see \cite  {Kh99b,He05}), we have $<y_s, \sigma  (x_s, x_{s+1}, \cdots, x_{\frac {k-s}{2}  } )
> \sigma  (x_{\frac {k-s}{2} +1}, x_{\frac {k-s}{2}+2}, \cdots, x_{k-1} ) \not=0$,
which implies $< y_s, \sigma  (x_s, x_{s+1}, \cdots, x_{k-1})> \not=0$.
\hfill $\Box$

\begin {Example} \label{4.2} {\rm (i)} Assume that $V_m$ is a  quantum vector space with $q_{i,i+1} \not=1$ and $q_{ij}=1$, $i+1 <j\le m, 1 \le i <m$. Let $A_k := \oplus _{i=1}^k \mathfrak B (V_i) $
and $A := \oplus _{i=1}^\infty \mathfrak B (V_i)$ as associative algebras. Then $A^-$ is not a  solvable Lie algebra.

 {\rm (ii)} $A^-$ has no  maximal nilpotent ideal and maximal solvable ideal.

\end {Example}

\noindent {\bf Proof.} {\rm (i) } If $A^-$ is solvable, then there
 exists a natural number $m$ such that $D^{(m)}(A^-)=0$. Therefore $D^{(m)}( \mathfrak B (V_{2^m}) ^- )=0$ and $ \sigma (x_1, x_2, \cdots, x_ {2^m}) \in D^{(m)}( \mathfrak B (V_{2^m})  ^-),$ which contradicts to  Lemma \ref {4.1}.

 {\rm (ii)} It follows from {\rm (i)}.
\hfill $\Box$
\vskip.1in
Similarly, we have the following conclusion.
 \begin {Example} \label{4.3} {\rm (i)}  Assume that $V_m$ is a  quantum vector space with $q_{i,i+1} \not=1$ and $q_{ij}=1$, $i+1 <j\le m, 1 \le i <m$.  Let $L_k := \oplus _{i=1}^k \mathfrak L^- (V_i) $
and $L := \oplus _{i=1}^\infty \mathfrak L ^-(V_i)$ as Lie algebras. Then $L$ is not any solvable Lie algebras.

 {\rm (ii)} $L$ has no  maximal nilpotent ideal and maximal solvable ideal.

\end {Example}

For a coagebra $C$, let $C^+ := \{ x\in C \mid \epsilon (x) =0\}$ (see \cite {Mo93,Ka95,Sw69a}).

 \begin {Example} \label {4.4}
   Assume that $V_m$ is a  quantum vector space with $q_{i,i+1} \not=1$ and $q_{ij}=1$, $i+1 <j\le m, 1 \le i <m$. Let $A_k := \oplus _{i=1}^k \mathfrak B (V_i) $
and $A := \oplus _{i=1}^\infty \mathfrak B (V_i)$ as associative algebras. Then

{\rm (i)}  $A^+ $ is not a nilpotent ideal of associative  algebra $A$.

 {\rm (ii)} $A$ has no maximal nilpotent ideal.

 {\rm (iii)} $A^+$ is a Baer radical of $A$, where Baer ideal is defined in \cite {Sz82}.

\end {Example}
\noindent {\bf Proof.} {\rm (i)} It follows Lemma \ref {4.1}.

{\rm (ii)} If $A$ has the maximal ideal $I, $ then  $I = A^+$ since every $\mathfrak B (V_i)^+$ is nilpotent for any $i \in \mathbb N. $ This contradicts to {\rm (i)}.

 {\rm (iii)} $A^+ = \sum _{k =1} (A_k )^+  $
 and $(A_1)^+ \subseteq (A_2)^+\subseteq \cdots (A_n)^+ \subseteq \cdots $. Since $(A_k)^+$ is  nilpotent ideal of $A$ for $k \in \mathbb N,$ $A^+ $ is  Baer radical of $A$ by \cite {Sz82}.  \hfill $\Box$

\section {Appendix }\label {s6}

\subsection {No $m$-infinity elements in Nichols algebras over finite cyclic groups} \label {s4}

In this subsection of the appendix we show that there  does not exist any $m$-infinity element in Nichols algebra $\mathfrak B(V) $ with arithmetic root system $\Delta (\mathfrak B(V))$  over finite cyclic groups.

\begin {Lemma} \label {40.1}
{\rm (i)} $<y_j, u ^{k}>=\sum \limits _{i=0}^{k - 1}p_{x_{j},u}^{ - i} u ^{i}<y_j, u > u ^{k - i - 1}$ for homogeneous element $u$.

{\rm (ii)} If $p_{11}^2 = p_{33}^2 = 1$, then   $$<y_1,[1,3]^{k}>=( - p_{13}^{ - 1})^{k}((p_{13}p_{31})^{k} - 1)\stackrel{k} { \overbrace{3131\cdots 1313}}.$$

{\rm (iii)} If $p_{33}^2 = p_{22}^2 = 1$ and $[2, 3]=0,$ then  $$<y_1,[[1,3],2]^{k}>=
\sum \limits _{i=0}^{k - 1}(p_{11}p_{12}p_{13}p_{21}p_{31})^{ - i}(p_{12}^{ - 1} - p_{21})(p_{13}^{ - 1} - p_{31})(p_{21}p_{31})^{k - 1}
\stackrel{k} { \overbrace{231231\cdots 312}}3.$$

{\rm (iv)} If $<y_i,   v > =0$ and $v$ is a homogeneous element,   then
$<y_{i},  [   v,  x_i] >=0.$

{\rm (v)} $<y_2,[[1,3],2]^{k}>=0$ and  $<y_3,[[1,3],2]^{k}>=0$.
\end {Lemma}
\noindent {\bf Proof.}
{\rm (i)}  It can be obtained by induction on $k$.

{\rm (ii)}
\begin{eqnarray}
<y_1,[1,3]^{k}>&=&\sum \limits _{i=0}^{k - 1}p_{11}^{ - i}p_{13}^{ - i}[1,3]^{i}<y_1,[1,3]>[1,3]^{k - i - 1}\nonumber\\
&=&\sum \limits _{i=0}^{k - 1}p_{11}^{ - i}p_{13}^{ - i}(31)^{i}(p_{13}^{ - 1} - p_{31})3( - p_{31})^{k - i - 1}(13)^{k - i - 1}\nonumber\\
&=&\sum \limits _{i=0}^{k - 1}p_{11}^{ - i}p_{13}^{ - i}( - p_{31})^{ - i}(p_{13}^{ - 1} - p_{31})(31)^{i}3( - p_{31})^{k - 1}(13)^{k - i - 1}\nonumber\\
&=&\sum \limits _{i=0}^{k - 1}p_{13}^{ - i}p_{31}^{ - i}(p_{13}^{ - 1} - p_{31})( - p_{31})^{k - 1}\stackrel{k} { \overbrace{3131\cdots 1313}}\nonumber\\
&=&( - p_{13}^{ - 1})^{k}((p_{13}p_{31})^{k} - 1)\stackrel{k} { \overbrace{3131\cdots 1313}}\nonumber.
\end{eqnarray}

{\rm (iii)}
\begin{eqnarray}
&&<y_1,[[1,3],2]^{k}>=\sum \limits _{i=0}^{k - 1}p_{11}^{ - i}p_{12}^{ - i}p_{13}^{ - i}[[1,3],2]^{i}<y_1,[[1,3],2]>[[1,3],2]^{k - i - 1}\nonumber\\
&&~~~~=\sum \limits _{i=0}^{k - 1}p_{11}^{ - i}p_{12}^{ - i}p_{13}^{ - i}[[1,3],2]^{i}(p_{12}^{ - 1} - p_{21})(p_{13}^{ - 1} - p_{31})23[[1,3],2]^{k - i - 1}\nonumber\\
&&~~~~=\sum \limits _{i=0}^{k - 1}p_{11}^{ - i}p_{12}^{ - i}p_{13}^{ - i}(231)^{i}(p_{12}^{ - 1} - p_{21})(p_{13}^{ - 1} - p_{31})23
(p_{21}p_{23}p_{31})^{k - i - 1}(132)^{k - i - 1}\nonumber\\
&&~~~~=\sum \limits _{i=0}^{k - 1}p_{11}^{ - i}p_{12}^{ - i}p_{13}^{ - i}(p_{21}p_{31})^{ - i}(p_{12}^{ - 1} - p_{21})(p_{13}^{ - 1} - p_{31})(p_{21}p_{31})^{k - 1}
(231)^{i}23(123)^{k - i - 1}\nonumber\\
&&~~~~=\sum \limits _{i=0}^{k - 1}(p_{11}p_{12}p_{13}p_{21}p_{31})^{ - i}(p_{12}^{ - 1} - p_{21})(p_{13}^{ - 1} - p_{31})(p_{21}p_{31})^{k - 1}
\stackrel{k} { \overbrace{231231\cdots 312}}3\nonumber.
\end{eqnarray}

{\rm (iv)}    It is clear.

{\rm (v) } It follows from {\rm (iv)}. \hfill $\Box$

\begin {Theorem} \label {40.2} If  $\Delta (\mathfrak B(V))$ is an arithmetic root system and
$V$ is  a ${\rm YD}$ module over finite cyclic group, then
 there does not exist any $m$-infinity element in  $\mathfrak B(V).$
\end {Theorem}

\noindent {\bf Proof.} We show this by several steps. Obviously,   it follows that
 $[x_1,  x_2]$,   $[x_1,  x_3]$ and $[[x_1,  x_3],  x_2]$ are nilpotent from  Lemma \ref {40.1} or simple computation. We only need consider the three cases below by Theorem 2.7 in \cite {WZZ14}.

{\rm (i)}  Assume generalized Dynkin diagram of braided vector space $V$ is the following condition:

  $\begin{picture}(100,    20) \put(27,   1){\makebox(0,
0)[t]{$\bullet$}} \put(60,    1){\makebox(0,   0)[t]{$\bullet$}}
\put(93,    1){\makebox(0,    0)[t]{$\bullet$}} \put(28,
 - 1){\line(1, 0){33}} \put(61,     - 1){\line(1,    0){30}} \put(18,
7){$ - 1$} \put(35,    6){$q$} \put(58, 7){$ - 1$}  \put(75,
6){$q^{ - 1}$}  \put(93,    7){$ - 1$} \put(110,    1){$, $}\ \ \
\put(120,    1) {$ q \in R_m, m>2.$}
\end{picture}$ \\
Then $D:=\{ [x_1]; [x_2]; [x_3]; [x_1,x_2]; [x_1,x_3];[[x_1,x_3],x_2]\}$ by Theorem 3.3 (i) in
\cite {WZZ14}. It follows that  every hard super-letter is nilpotent from Lemma
\ref {40.1}.

{\rm (ii)}   Assume generalized Dynkin diagram of braided vector space $V$ is the following condition:

 $\begin{picture}(100,    20) \put(27,   1){\makebox(0,
0)[t]{$\bullet$}} \put(60,    1){\makebox(0,   0)[t]{$\bullet$}}
\put(93,    1){\makebox(0,    0)[t]{$\bullet$}} \put(28,
 - 1){\line(1, 0){33}} \put(61,     - 1){\line(1,    0){30}} \put(18,
7){$ - 1$} \put(35,    6){$\zeta $} \put(58, 7){$ - 1$}  \put(75,
6){$\zeta$}  \put(93,    7){${ - 1}$} \put(110,    1){$, $}\ \ \
\put(120,    1) {$ \zeta \in R_3.$}
\end{picture}$ \\
Then $D =\{ [x_1]; [x_2]; [x_3]; [x_1,x_2]; [x_1,x_3];
[[x_1,x_3],x_2]; [[x_1,x_2],[x_1,x_3]]; [[x_1,x_2],[[x_1,x_3],x_2]]$;
$[[x_1,x_3],[[x_1,x_3],x_2]]; [[[x_1,x_2],[[x_1,x_3],x_2]],[x_1,x_3]]\}$  by Theorem 3.3(i) in \cite  {WZZ14}.

\vskip.1in
Now we show  that every hard super-letter is nilpotent.
\begin{eqnarray}
&&<y_1,[[1,2],[1,3]]>= - p_{12}^{ - 1}p_{13}^{ - 1}(1 - p_{12}p_{21}p_{31}p_{13})[2,[1,3]],\nonumber\\
&&<y_1,[[1,2],[1,3]]^{3}>=\sum \limits _{i=0}^{2}p_{11}^{ - 2i}p_{12}^{ - i}p_{13}^{ - i}[[1,2],[1,3]]^{i}<y_1,[[1,2],[1,3]]>[[1,2],[1,3]]^{2 - i}\nonumber\\
&&~~~~=\sum \limits _{i=0}^{2}p_{12}^{ - i}p_{13}^{ - i}[[1,2],[1,3]]^{i}( - p_{12}^{ - 1}p_{13}^{ - 1}
(1 - p_{12}p_{21}p_{31}p_{13})[2,[1,3]])[[1,2],[1,3]]^{2 - i}\nonumber\\
&&~~~~= - p_{12}^{ - 1}p_{13}^{ - 1}(1 - p_{12}p_{21}p_{31}p_{13})\{[2,[1,3]][[1,2],[1,3]]^{2} \nonumber\\
&&~~~~~+  p_{12}^{ - 1}p_{13}^{ - 1}[[1,2],[1,3]][2,[1,3]]
   [[1,2],[1,3]] +  p_{12}^{ - 2}p_{13}^{ - 2}[[1,2],[1,3]]^{2}[2,[1,3]]\}\nonumber\\
&&~~~~:= - p_{12}^{ - 1}p_{13}^{ - 1}(1 - p_{12}p_{21}p_{31}p_{13})R,\nonumber
\end{eqnarray}
\begin{eqnarray}
&&<y_2,[2,[1,3]]>=p_{21}^{ - 1}p_{23}^{ - 1}(1 - p_{12}p_{21})[1,3],\nonumber\\
&&<y_2,[2,[1,3]]^2>=p_{21}^{ - 1}p_{23}^{ - 1}(1 - p_{12}p_{21})\{[1,3][2,[1,3]]
+  p_{21}^{ - 1}p_{22}^{ - 1}p_{23}^{ - 1}[2,[1,3]][1,3]\},\nonumber
\end{eqnarray}
\begin{eqnarray}
&&<y_2,R>=p_{21}^{ - 1}p_{23}^{ - 1}(1 - p_{12}p_{21})\{[1,3][[1,2],[1,3]]^{2} \nonumber\\
&&~~~~~+  p_{12}^{ - 1}p_{13}^{ - 1}p_{21}^{ - 2}
p_{22}^{ - 1}p_{23}^{ - 1}[[1,2],[1,3]][1,3][[1,2],[1,3]]
 +  p_{21}^{ - 4}p_{22}^{ - 2}p_{23}^{ - 2}p_{12}^{ - 2}p_{13}^{ - 2}[[1,2],[1,3]]^{2}[1,3]\}\nonumber\\
&&~~~~=p_{21}^{ - 1}p_{23}^{ - 1}(1 - p_{12}p_{21})
\{p_{11}^{4}p_{12}^{2}p_{13}^{2}p_{31}^{4}p_{32}^{2}p_{33}^{2}[[1,2],[1,3]]^{2}[1,3]\nonumber\\
&&~~~~~ +  p_{11}^{2}p_{12}p_{13}p_{31}^{2}p_{32}p_{33}p_{12}^{ - 1}p_{13}^{ - 1}p_{21}^{ - 2}
p_{22}^{ - 1}p_{23}^{ - 1}[[1,2],[1,3]][[1,2],[1,3]][1,3]\nonumber\\
&&~~~~~ +  p_{21}^{ - 4}p_{22}^{ - 2}p_{23}^{ - 2}p_{12}^{ - 2}p_{13}^{ - 2}[[1,2],[1,3]]^{2}[1,3]\}\nonumber\\
&&~~~~=p_{21}^{ - 1}p_{23}^{ - 1}(1 - p_{12}p_{21})p_{21}^{ - 2}p_{32}^{2}p_{31}^{2}
\{p_{12}^{2}p_{21}^{2}p_{13}^{2}p_{31}^{2}
 +  1 +  p_{21}^{ - 2}p_{13}^{ - 2}p_{12}^{ - 2}p_{31}^{ - 2}\}[[1,2],[1,3]]^{2}[1,3]=0\nonumber\\
&&~~~~~~~~{\rm since}~ [[[1,2],[1,3]],[1,3]]=0,\nonumber
\end{eqnarray}
\begin{eqnarray}
&& <y_3,R>=0,\nonumber\\
&&<y_1,R>=<y_1,<y_1,[[1,2],[1,3]]^{3}>>=<y_{1}^2,[[1,2],[1,3]]^{3}>=0. {\rm Then}~ [[1,2],[1,3]]^{3}=0,\nonumber\\
&&<y_1,[[1,2],[[1,3],2]]>= - p_{12}^{ - 1}p_{21}p_{13}^{ - 1}(1 - \xi)^{2}2[1,3]2,\nonumber
\end{eqnarray}
\begin{eqnarray}
&&<y_1,[[1,2],[[1,3],2]]^2>= - p_{12}^{ - 1}p_{21}p_{13}^{ - 1}(1 - \xi)^{2}2[1,3]2[[1,2],[[1,3],2]]\nonumber\\
&&~~~~~ +  p_{11}^{ - 2}p_{12}^{ - 2}p_{13}^{ - 1}[[1,2],[[1,3],2]]( - p_{12}^{ - 1}p_{21}p_{13}^{ - 1}
(1 - \xi)^{2}2[1,3]2)\nonumber\\
&&~~~~= - p_{12}^{ - 1}p_{21}p_{13}^{ - 1}(1 - \xi)^{2}\cdot\{2[1,3]2[[1,2],[[1,3],2]]
 +  p_{11}^{ - 2}p_{12}^{ - 2}p_{13}^{ - 1}[[1,2],[[1,3],2]]2[1,3]2\}\nonumber\\
&&~~~~:= - p_{12}^{ - 1}p_{21}p_{13}^{ - 1}(1 - \xi)^{2}\cdot A.\nonumber
\end{eqnarray}
It is clear $<y_1,A>=0,<y_3,A>=0$, and
\begin{eqnarray}
&&<y_2,A>=[1,3]2[[1,2],[[1,3],2]] - p_{21}^{ - 1}p_{23}^{ - 1}2[1,3][[1,2],[[1,3],2]]\nonumber\\
&&~~~~~ +  \xi p_{13}^{ - 1}p_{23}^{ - 1}[[1,2],[[1,3],2]][1,3]2
 - \xi p_{21}^{ - 1}p_{13}^{ - 1}p_{23}^{ - 2}[[1,2],[[1,3],2]]2[1,3]\nonumber\\
&&~~~~= - p_{21}^{ - 1}p_{23}^{ - 1}[[1,3],2][[1,2],[[1,3],2]]
 - \xi p_{21}^{ - 1}p_{13}^{ - 1}p_{23}^{ - 2}[[1,2],[[1,3],2]][[1,3],2]\nonumber\\
&&~~~~= - p_{21}^{ - 1}p_{23}^{ - 1}([[1,3],2][[1,2],[[1,3],2]]
 +  \xi p_{13}^{ - 1}p_{23}^{ - 1}[[1,2],[[1,3],2]][[1,3],2])\nonumber\\
&&~~~~= - p_{21}^{ - 1}p_{23}^{ - 1}[[[1,2],[[1,3],2]],[[1,3],2]]=0, \nonumber
\end{eqnarray}
then $[[1,2],[[1,3],2]]^2=0$ and 
$<y_1,[[1,3],[[1,3],2]]>= - p_{13}^{ - 2}p_{21}p_{23}(1 - \xi)^{2}3[1,2]3$.

We obtain $[[1,3],[[1,3],2]]^2=0$ since $[[1,3],2]=q_{32}^{ - 1}[[1,2],3]$.
\begin{eqnarray}
&&<y_1,[[[1,2],[[1,3],2]],[1,3]]>=p_{12}^{ - 2}p_{13}^{ - 2}\xi(1 - \xi)[2,[1,3]]^{2},\nonumber\\
&&<y_1,[[[1,2],[[1,3],2]],[1,3]]^2>=p_{12}^{ - 2}p_{13}^{ - 2}\xi(1 - \xi)[2,[1,3]]^{2}[[[1,2],[[1,3],2]],[1,3]]\nonumber\\
&&~~~~~ +  p_{11}^{ - 3}p_{12}^{ - 2}p_{13}^{ - 2}[[[1,2],[[1,3],2]],[1,3]](p_{12}^{ - 2}p_{13}^{ - 2}
\xi(1 - \xi)[2,[1,3]]^{2})\nonumber\\
&&~~~~=p_{12}^{ - 2}p_{13}^{ - 2}\xi(1 - \xi)\{[2,[1,3]]^{2}[[[1,2],[[1,3],2]],[1,3]]\nonumber\\
&&~~~~~~~ +  p_{11}^{ - 3}p_{12}^{ - 2}p_{13}^{ - 2}[[[1,2],[[1,3],2]],[1,3]][2,[1,3]]^{2})\}\nonumber\\
&&~~~~:=p_{12}^{ - 2}p_{13}^{ - 2}\xi(1 - \xi)B,\nonumber
\end{eqnarray}
\begin{eqnarray}
&&<y_2,B>=p_{21}^{ - 1}p_{23}^{ - 1}(1 - p_{12}p_{21})([1,3][2,[1,3]]
  +  p_{21}^{ - 1}p_{22}^{ - 1}p_{23}^{ - 1}[2,[1,3]][1,3])[[[1,2],[[1,3],2]],[1,3]]\nonumber\\
&&~~~~~ +  p_{21}^{ - 3}p_{22}^{ - 2}p_{23}^{ - 2}p_{11}^{ - 3}p_{12}^{ - 2}p_{13}^{ - 2}[[[1,2],[[1,3],2]],[1,3]]p_{21}^{ - 1}p_{23}^{ - 1}(1 - p_{12}p_{21})\nonumber\\
&&~~~~~~~~~([1,3][2,[1,3]] +  p_{21}^{ - 1}p_{22}^{ - 1}p_{23}^{ - 1}[2,[1,3]][1,3])\nonumber\\
&&~~~~~~~~~:=C,\nonumber
\end{eqnarray}
\begin{eqnarray}
&&<y_1,C>=p_{21}^{ - 1}p_{23}^{ - 1}(1 - p_{12}p_{21})((p_{13}^{ - 1} - p_{31})3[2,[1,3]] \nonumber\\
&&~~~~~+  p_{11}^{ - 1}p_{12}^{ - 1}p_{13}^{ - 1}
p_{21}^{ - 1}p_{22}^{ - 1}p_{23}^{ - 1}[2,[1,3]](p_{13}^{ - 1} - p_{31})3)[[[1,2],[[1,3],2]],[1,3]]\nonumber\\
&&~~~~~ +  p_{12}^{ - 1}p_{13}^{ - 2}p_{21}^{ - 1}p_{23}^{ - 1}(1 - p_{12}p_{21})([1,3][2,[1,3]]\nonumber\\
&&~~~~~~~+  p_{21}^{ - 1}p_{22}^{ - 1}p_{23}^{ - 1}[2,[1,3]][1,3])
p_{12}^{ - 2}p_{13}^{ - 2}\xi(1 - \xi)[2,[1,3]]^{2}\nonumber\\
&&~~~~~ +  p_{21}^{ - 3}p_{22}^{ - 2}p_{23}^{ - 2}p_{11}^{ - 3}p_{12}^{ - 2}p_{13}^{ - 2}
p_{12}^{ - 2}p_{13}^{ - 2}\xi(1 - \xi)[2,[1,3]]^{2}
p_{21}^{ - 1}p_{23}^{ - 1}(1 - p_{12}p_{21})([1,3][2,[1,3]]\nonumber\\
&&~~~~~ +  p_{21}^{ - 1}p_{22}^{ - 1}p_{23}^{ - 1}[2,[1,3]][1,3]) +  p_{21}^{ - 3}p_{22}^{ - 2}
p_{23}^{ - 2}p_{12}^{ - 4}p_{13}^{ - 4}
[[[1,2],[[1,3],2]],[1,3]]p_{21}^{ - 1}p_{23}^{ - 1}\nonumber\\
&&~~~~~~(1 - p_{12}p_{21})
((p_{13}^{ - 1} - p_{31})3[2,[1,3]] +  p_{11}^{ - 1}p_{12}^{ - 1}p_{13}^{ - 1}p_{21}^{ - 1}p_{22}^{ - 1}p_{23}^{ - 1}[2,[1,3]](p_{13}^{ - 1} - p_{31})3)\nonumber\\
&&~~~~=p_{21}^{ - 1}p_{23}^{ - 1}p_{13}^{ - 1}(1 - \xi)^{2}(3[2,[1,3]] +  p_{12}^{ - 1}p_{13}^{ - 1}
p_{21}^{ - 1}p_{23}^{ - 1}[2,[1,3]]3)[[[1,2],[[1,3],2]],[1,3]]\nonumber\\
&&~~~~~ +  p_{12}^{ - 3}p_{13}^{ - 4}p_{21}^{ - 1}p_{23}^{ - 1}\xi(1 - \xi)^{2}
([1,3][2,[1,3]] - p_{21}^{ - 1}p_{23}^{ - 1}[2,[1,3]][1,3])[2,[1,3]]^{2}\nonumber\\
&&~~~~~ - p_{13}^{ - 4}p_{23}^{ - 3}(1 - \xi)^{2}[2,[1,3]]^{2}
([1,3][2,[1,3]] - p_{21}^{ - 1}p_{23}^{ - 1}[2,[1,3]][1,3])\nonumber\\
&&~~~~~ +  p_{23}^{ - 3}p_{12}^{ - 1}p_{13}^{ - 5}p_{21}^{ - 1}(1 - \xi)^{2}
[[[1,2],[[1,3],2]],[1,3]](3[2,[1,3]] +  p_{12}^{ - 1}p_{13}^{ - 1}p_{21}^{ - 1}p_{23}^{ - 1}[2,[1,3]]3)\nonumber\\
&&~~~~:=D,\nonumber
\end{eqnarray}
\begin{eqnarray}
&&<y_3,D>=p_{21}^{ - 1}p_{23}^{ - 1}p_{13}^{ - 1}(1 - \xi)^{2}([2,[1,3]]\nonumber\\
&&~~~~~~~~ +  p_{12}^{ - 1}p_{13}^{ - 1}
p_{21}^{ - 1}p_{23}^{ - 1}p_{31}^{ - 1}p_{32}^{ - 1}p_{33}^{ - 1}[2,[1,3]])[[[1,2],[[1,3],2]],[1,3]]\nonumber\\
&&~~~~~ +  p_{23}^{ - 3}p_{12}^{ - 1}p_{13}^{ - 5}p_{21}^{ - 1}(1 - \xi)^{2}p_{31}^{ - 3}p_{32}^{ - 2}
[[[1,2],[[1,3],2]],[1,3]]([2,[1,3]] \nonumber\\
&&~~~~~~~+  p_{12}^{ - 1}p_{13}^{ - 1}
p_{21}^{ - 1}p_{23}^{ - 1}p_{31}^{ - 1}p_{32}^{ - 1}p_{33}^{ - 1}[2,[1,3]])\nonumber\\
&&~~~~=p_{21}^{ - 1}p_{23}^{ - 1}p_{13}^{ - 1}(1 - \xi)^{3}[2,[1,3]][[[1,2],[[1,3],2]],[1,3]]\nonumber\\
&&~~~~~ +  p_{23}^{ - 1}p_{12}^{ - 1}p_{13}^{ - 2}p_{21}^{ - 1}(1 - \xi)^{3}[[[1,2],[[1,3],2]],[1,3]][2,[1,3]]\nonumber\\
&&~~~~:=E,\nonumber
\end{eqnarray}
\begin{eqnarray}
&&<y_1,E>=p_{21}^{ - 1}p_{23}^{ - 1}p_{13}^{ - 1}p_{11}^{ - 1}p_{12}^{ - 1}
p_{13}^{ - 1}(1 - \xi)^{3}[2,[1,3]]
p_{12}^{ - 2}p_{13}^{ - 2}\xi(1 - \xi)[2,[1,3]]^{2}\nonumber\\
&&~~~~~ +  p_{23}^{ - 1}p_{12}^{ - 1}p_{13}^{ - 2}p_{21}^{ - 1}(1 - \xi)^{3}p_{12}^{ - 2}p_{13}^{ - 2}
\xi(1 - \xi)[2,[1,3]]^{2}[2,[1,3]]\nonumber\\
&&~~~~= - p_{23}^{ - 1}p_{13}^{ - 4}(1 - \xi)^{4}p_{12}^{ - 2}[2,[1,3]]^{3}
+  p_{23}^{ - 1}p_{13}^{ - 4}(1 - \xi)^{4}p_{12}^{ - 2}[2,[1,3]]^{3}=0,   \nonumber
\end{eqnarray}
\begin{eqnarray}
&&<y_2,E>=p_{21}^{ - 1}p_{23}^{ - 1}p_{13}^{ - 1}(1 - \xi)^{3}p_{21}^{ - 1}p_{23}^{ - 1}(1 - p_{12}p_{21})[1,3][[[1,2],[[1,3],2]],[1,3]]\nonumber\\
&&~~~~~ +  p_{21}^{ - 3}p_{23}^{ - 2}p_{23}^{ - 1}p_{12}^{ - 1}p_{13}^{ - 2}p_{21}^{ - 1}(1 - \xi)^{3}p_{21}^{ - 1}p_{23}^{ - 1}(1 - p_{12}p_{21})
[[[1,2],[[1,3],2]],[1,3]][1,3]\nonumber\\
&&~~~~=p_{21}^{ - 2}p_{23}^{ - 2}p_{13}^{ - 1}(1 - \xi)^{4}[1,3][[[1,2],[[1,3],2]],[1,3]]\nonumber\\
&&~~~~~~~ +  p_{21}^{ - 5}p_{23}^{ - 4}p_{12}^{ - 1}p_{13}^{ - 2}
(1 - \xi)^{4}[[[1,2],[[1,3],2]],[1,3]][1,3]\nonumber\\
&&~~~~=p_{21}^{ - 2}p_{23}^{ - 2}p_{13}^{ - 1}(1 - \xi)^{4}\{[1,3][[[1,2],[[1,3],2]],[1,3]]\nonumber\\
&&~~~~~~~ +  p_{21}^{ - 3}p_{23}^{ - 2}p_{12}^{ - 1}p_{13}^{ - 1}[[[1,2],[[1,3],2]],[1,3]][1,3]\}\nonumber\\
&&~~~~=0 ~{\rm  since}~[[[[1,2],[[1,3],2]],[1,3]],[1,3]]=0,\nonumber
\end{eqnarray}
\begin{eqnarray}
&&<y_2,D>=p_{21}^{ - 1}p_{23}^{ - 1}p_{13}^{ - 1}(1 - \xi)^{2}(p_{23}^{ - 1}3p_{21}^{ - 1}p_{23}^{ - 1}(1 - p_{12}p_{21})[1,3] \nonumber\\
&&~~~~~+  p_{12}^{ - 1}p_{13}^{ - 1}
p_{21}^{ - 1}p_{23}^{ - 1}p_{21}^{ - 1}p_{23}^{ - 1}(1 - p_{12}p_{21})[1,3]3)[[[1,2],[[1,3],2]],[1,3]]\nonumber\\
&&~~~~~ - p_{12}^{ - 3}p_{13}^{ - 4}p_{21}^{ - 1}p_{23}^{ - 1}\xi(1 - \xi)^{2}p_{21}^{ - 2}p_{23}^{ - 2}
([1,3][2,[1,3]]\nonumber\\
&&~~~~~ - p_{21}^{ - 1}p_{23}^{ - 1}[2,[1,3]][1,3])p_{21}^{ - 1}p_{23}^{ - 1}(1 - p_{12}p_{21})
([1,3][2,[1,3]] +  p_{21}^{ - 1}p_{22}^{ - 1}p_{23}^{ - 1}[2,[1,3]][1,3])\nonumber\\
&&~~~~~ - p_{13}^{ - 4}p_{23}^{ - 3}(1 - \xi)^{2}p_{21}^{ - 1}
p_{23}^{ - 1}(1 - p_{12}p_{21})([1,3][2,[1,3]] +  p_{21}^{ - 1}p_{22}^{ - 1}p_{23}^{ - 1}[2,[1,3]][1,3])\nonumber\\
&&~~~~~~~~([1,3][2,[1,3]] - p_{21}^{ - 1}p_{23}^{ - 1}[2,[1,3]][1,3])\nonumber\\
&&~~~~~ +  p_{21}^{ - 3}p_{23}^{ - 2}p_{23}^{ - 3}p_{12}^{ - 1}p_{13}^{ - 5}p_{21}^{ - 1}(1 - \xi)^{2}
[[[1,2],[[1,3],2]],[1,3]](p_{23}^{ - 1}3p_{21}^{ - 1}p_{23}^{ - 1}(1 - p_{12}p_{21})[1,3]\nonumber\\
&&~~~~~ +  p_{12}^{ - 1}p_{13}^{ - 1}
p_{21}^{ - 1}p_{23}^{ - 1}p_{21}^{ - 1}p_{23}^{ - 1}(1 - p_{12}p_{21})[1,3]3)\nonumber\\
&&~~~~=p_{21}^{ - 2}p_{23}^{ - 3}p_{13}^{ - 1}(1 - \xi)^{3}(3[1,3] +  p_{12}^{ - 1}p_{13}^{ - 1}
p_{21}^{ - 1}[1,3]3)[[[1,2],[[1,3],2]],[1,3]] \nonumber\\
&&~~~~~~~- p_{12}^{ - 3}p_{13}^{ - 4}p_{21}^{ - 4}p_{23}^{ - 4}\xi(1 - \xi)^{3}
([1,3][2,[1,3]] \nonumber\\
&&~~~~~~~- p_{21}^{ - 1}p_{23}^{ - 1}[2,[1,3]][1,3])
([1,3][2,[1,3]] +  p_{21}^{ - 1}p_{22}^{ - 1}p_{23}^{ - 1}[2,[1,3]][1,3]) \nonumber\\
&&~~~~~- p_{13}^{ - 4}p_{23}^{ - 4}(1 - \xi)^{3}p_{21}^{ - 1}
([1,3][2,[1,3]] +  p_{21}^{ - 1}p_{22}^{ - 1}p_{23}^{ - 1}[2,[1,3]][1,3])([1,3][2,[1,3]] \nonumber\\
&&~~~~~~~~~~- p_{21}^{ - 1}p_{23}^{ - 1}[2,[1,3]][1,3])\nonumber\\
&&~~~~~ +  p_{23}^{ - 7}p_{12}^{ - 1}p_{13}^{ - 5}p_{21}^{ - 5}(1 - \xi)^{3}
[[[1,2],[[1,3],2]],[1,3]](3[1,3] +  p_{12}^{ - 1}p_{13}^{ - 1}p_{21}^{ - 1}[1,3]3)\nonumber
\end{eqnarray}
\begin{eqnarray}
&&~~~=~p_{21}^{ - 2}p_{23}^{ - 3}p_{13}^{ - 1}(1 - \xi)^{4}3[1,3][[[1,2],[[1,3],2]],[1,3]]\nonumber\\
&&~~~~~ - p_{13}^{ - 4}p_{21}^{ - 1}p_{23}^{ - 4}\xi(1 - \xi)^{3}
([1,3][2,[1,3]] - p_{21}^{ - 1}p_{23}^{ - 1}[2,[1,3]][1,3])
([1,3][2,[1,3]]\nonumber\\
&&~~~~~~~~~~ - p_{21}^{ - 1}p_{23}^{ - 1}[2,[1,3]][1,3])\nonumber\\
&&~~~~~- p_{13}^{ - 4}p_{23}^{ - 4}(1 - \xi)^{3}p_{21}^{ - 1}
([1,3][2,[1,3]] - p_{21}^{ - 1}p_{23}^{ - 1}[2,[1,3]][1,3])
([1,3][2,[1,3]]\nonumber\\
&&~~~~~~~~~~ - p_{21}^{ - 1}p_{23}^{ - 1}[2,[1,3]][1,3])\nonumber\\
&&~~~~~ +  p_{23}^{ - 7}p_{12}^{ - 1}p_{13}^{ - 5}p_{21}^{ - 5} (1 - \xi)^{4}
[[[1,2],[[1,3],2]],[1,3]]3[1,3]\nonumber
\end{eqnarray}
\begin{eqnarray}
&&~~~=~p_{21}^{ - 2}p_{23}^{ - 3}p_{13}^{ - 1} (1 - \xi)^{4}3[1,3][[[1,2],[[1,3],2]],[1,3]]\nonumber\\
&&~~~~~+  p_{13}^{ - 4}p_{21}^{ - 1}p_{23}^{ - 4}\xi^{2} (1 - \xi)^{3}
 ([1,3][2,[1,3]] - p_{21}^{ - 1}p_{23}^{ - 1}[2,[1,3]][1,3])
 ([1,3][2,[1,3]]\nonumber\\
&&~~~~~~~~~~- p_{21}^{ - 1}p_{23}^{ - 1}[2,[1,3]][1,3])\nonumber\\
&&~~~~~+  p_{23}^{ - 7}p_{12}^{ - 1}p_{13}^{ - 5}p_{21}^{ - 5} (1 - \xi)^{4}
[[[1,2],[[1,3],2]],[1,3]]3[1,3]\nonumber\\
&&~~~:=~F,\nonumber
\end{eqnarray}
\begin{eqnarray}
&&<y_3,F>=p_{21}^{ - 2}p_{23}^{ - 3}p_{13}^{ - 1} (1 - \xi)^{4}[1,3][[[1,2],[[1,3],2]],[1,3]]\nonumber\\
&&~~~~~+  p_{31}^{ - 3}p_{32}^{ - 2}p_{23}^{ - 7}p_{12}^{ - 1}p_{13}^{ - 5}p_{21}^{ - 5} (1 - \xi)^{4}
[[[1,2],[[1,3],2]],[1,3]][1,3]\nonumber\\
&&~~~~=p_{21}^{ - 2}p_{23}^{ - 3}p_{13}^{ - 1} (1 - \xi)^{4} ([1,3][[[1,2],[[1,3],2]],[1,3]]\nonumber\\
&&~~~~~~ +  p_{12}^{ - 1}p_{13}^{ - 1}p_{21}^{ - 3}p_{23}^{ - 2}[[[1,2],[[1,3],2]],[1,3]][1,3])\nonumber\\
&&~~~~=0 ~{\rm since} ~[[[[1,2],[[1,3],2]],[1,3]],[1,3]]=0,\nonumber
\end{eqnarray}
\begin{eqnarray}
&&<y_1,F>= - p_{21}^{ - 2}p_{23}^{ - 3}p_{13}^{ - 3} (1 - \xi)^{4}3[1,3]p_{12}^{ - 2}p_{13}^{ - 2}\xi (1 - \xi)[2,[1,3]]^{2}\nonumber\\
&&~~~~~+  p_{13}^{ - 4}p_{21}^{ - 1}p_{23}^{ - 4}\xi^{2} (1 - \xi)^{3}
 ( (p_{13}^{ - 1} - p_{31})3[2,[1,3]] +  p_{21}^{ - 1}p_{23}^{ - 1}p_{12}^{ - 1}p_{13}^{ - 1}[2,[1,3]] (p_{13}^{ - 1} - p_{31})3)\nonumber\\
&&~~~~~~~~ ([1,3][2,[1,3]] - p_{21}^{ - 1}p_{23}^{ - 1}[2,[1,3]][1,3]) \nonumber\\
&&~~~~~+  p_{12}^{ - 1}p_{13}^{ - 2}p_{13}^{ - 4}p_{21}^{ - 1}p_{23}^{ - 4}\xi^{2} (1 - \xi)^{3}
 ([1,3][2,[1,3]] - p_{21}^{ - 1}p_{23}^{ - 1}[2,[1,3]][1,3])\nonumber\\
&&~~~~~~~~ ( (p_{13}^{ - 1} - p_{31})3[2,[1,3]] +  p_{21}^{ - 1}p_{23}^{ - 1}p_{12}^{ - 1}p_{13}^{ - 1}[2,[1,3]] (p_{13}^{ - 1} - p_{31})3)\nonumber\\
&&~~~~~ +  p_{23}^{ - 7}p_{12}^{ - 1}p_{13}^{ - 5}p_{21}^{ - 5} (1 - \xi)^{4}
p_{12}^{ - 2}p_{13}^{ - 2}\xi (1 - \xi)[2,[1,3]]^{2}3[1,3]\nonumber\\
&&~~~~= - p_{23}^{ - 3}p_{13}^{ - 5} (1 - \xi)^{5}\xi^{2}3[1,3][2,[1,3]]^{2}\nonumber\\
&&~~~~~ +  p_{13}^{ - 5}p_{21}^{ - 1}p_{23}^{ - 4}\xi^{2} (1 - \xi)^{4}
 (3[2,[1,3]] +  p_{21}^{ - 1}p_{23}^{ - 1}p_{12}^{ - 1}p_{13}^{ - 1}[2,[1,3]]3) \nonumber\\
&&~~~~~~~([1,3][2,[1,3]] - p_{21}^{ - 1}p_{23}^{ - 1}[2,[1,3]][1,3])\nonumber\\
&&~~~~~ +  p_{12}^{ - 1}p_{13}^{ - 7}p_{21}^{ - 1}p_{23}^{ - 4}\xi^{2} (1 - \xi)^{4}
 ([1,3][2,[1,3]] - p_{21}^{ - 1}p_{23}^{ - 1}[2,[1,3]][1,3])
 (3[2,[1,3]] \nonumber\\
&&~~~~~~~~+  p_{21}^{ - 1}p_{23}^{ - 1}p_{12}^{ - 1}p_{13}^{ - 1}[2,[1,3]]3)\nonumber\\
&&~~~~~ +  p_{23}^{ - 7}p_{13}^{ - 7}p_{21}^{ - 2} (1 - \xi)^{5}\xi[2,[1,3]]^{2}3[1,3]\nonumber\\
&&~~~~:=G,\nonumber
\end{eqnarray}
\begin{eqnarray}
&&<y_3,G>= - p_{23}^{ - 3}p_{13}^{ - 5} (1 - \xi)^{5}\xi^{2}[1,3][2,[1,3]]^{2}\nonumber\\
&&~~~~~ +  p_{13}^{ - 5}p_{21}^{ - 1}p_{23}^{ - 4}\xi^{2} (1 - \xi)^{4}
 ([2,[1,3]] - p_{31}^{ - 1}p_{32}^{ - 1}p_{21}^{ - 1}p_{23}^{ - 1}p_{12}^{ - 1}p_{13}^{ - 1}[2,[1,3]])\nonumber\\
&&~~~~~~~~ ([1,3][2,[1,3]] - p_{21}^{ - 1}p_{23}^{ - 1}[2,[1,3]][1,3])\nonumber\\
&&~~~~~ +  p_{12}^{ - 1}p_{13}^{ - 7}p_{21}^{ - 1}p_{23}^{ - 4}\xi^{2}
(1 - \xi)^{4}p_{31}^{ - 2}p_{32}^{ - 1}
 ([1,3][2,[1,3]] - p_{21}^{ - 1}p_{23}^{ - 1}[2,[1,3]][1,3])\nonumber\\
&&~~~~~~~~ ([2,[1,3]] - p_{31}^{ - 1}p_{32}^{ - 1}p_{21}^{ - 1}p_{23}^{ - 1}p_{12}^{ - 1}p_{13}^{ - 1}[2,[1,3]])\nonumber\\
&&~~~~~ +  p_{31}^{ - 2}p_{32}^{ - 2}p_{23}^{ - 7}p_{13}^{ - 7}p_{21}^{ - 2} (1 - \xi)^{5}\xi[2,[1,3]]^{2}[1,3]\nonumber\\
&&~~~~= - p_{23}^{ - 3}p_{13}^{ - 5} (1 - \xi)^{5}\xi^{2}[1,3][2,[1,3]]^{2}\nonumber\\
&&~~~~~ +  p_{13}^{ - 5}p_{21}^{ - 1}p_{23}^{ - 4}\xi^{2} (1 - \xi)^{5}[2,[1,3]] ([1,3][2,[1,3]] - p_{21}^{ - 1}p_{23}^{ - 1}[2,[1,3]][1,3])\nonumber\\
&&~~~~~ +  p_{13}^{ - 5}p_{23}^{ - 3}\xi^{2} (1 - \xi)^{5} ([1,3][2,[1,3]] - p_{21}^{ - 1}p_{23}^{ - 1}[2,[1,3]][1,3])[2,[1,3]]\nonumber\\
&&~~~~~ +  p_{23}^{ - 5}p_{13}^{ - 5}p_{21}^{ - 2} (1 - \xi)^{5}\xi^{2}[2,[1,3]]^{2}[1,3]\nonumber\\
&&~~~~= - p_{23}^{ - 3}p_{13}^{ - 5} (1 - \xi)^{5}\xi^{2}[1,3][2,[1,3]]^{2} +  p_{13}^{ - 5}p_{21}^{ - 1}p_{23}^{ - 4}\xi^{2} (1 - \xi)^{5}[2,[1,3]][1,3][2,[1,3]]\nonumber\\
&&~~~~~ - p_{13}^{ - 5}p_{21}^{ - 2}p_{23}^{ - 5}\xi^{2} (1 - \xi)^{5}[2,[1,3]]^{2}[1,3] +  p_{13}^{ - 5}p_{23}^{ - 3}\xi^{2} (1 - \xi)^{5}[1,3][2,[1,3]]^{2}\nonumber\\
&&~~~~~ - p_{13}^{ - 5}p_{23}^{ - 4}\xi^{2} (1 - \xi)^{5}p_{21}^{ - 1}[2,[1,3]][1,3][2,[1,3]]
 +  p_{23}^{ - 5}p_{13}^{ - 5}p_{21}^{ - 2} (1 - \xi)^{5}\xi^{2}[2,[1,3]]^{2}[1,3]=0,\nonumber
\end{eqnarray}
\begin{eqnarray}
&&<y_2,G>= - p_{23}^{ - 3}p_{13}^{ - 5} (1 - \xi)^{5}\xi^{2}p_{21}^{ - 1}p_{23}^{ - 2}3[1,3]p_{21}^{ - 1}p_{23}^{ - 1} (1 - p_{12}p_{21}) \nonumber\\
&&~~~~~~~~~([1,3][2,[1,3]]
 +  p_{21}^{ - 1}p_{22}^{ - 1}p_{23}^{ - 1}[2,[1,3]][1,3])\nonumber\\
&&~~~~~ +  p_{13}^{ - 5}p_{21}^{ - 1}p_{23}^{ - 4}\xi^{2} (1 - \xi)^{4}
 (p_{23}^{ - 1}3p_{21}^{ - 1}p_{23}^{ - 1} (1 - p_{12}p_{21})[1,3]\nonumber\\
&&~~~~~ +  p_{21}^{ - 1}p_{23}^{ - 1}p_{12}^{ - 1}p_{13}^{ - 1}p_{21}^{ - 1}p_{23}^{ - 1}
 (1 - p_{12}p_{21})[1,3]3) ([1,3][2,[1,3]] - p_{21}^{ - 1}p_{23}^{ - 1}[2,[1,3]][1,3])\nonumber\\
&&~~~~~ - p_{12}^{ - 1}p_{13}^{ - 7}p_{21}^{ - 1}p_{23}^{ - 4}\xi^{2} (1 - \xi)^{4}p_{21}^{ - 2}p_{23}^{ - 2}
 ([1,3][2,[1,3]] - p_{21}^{ - 1}p_{23}^{ - 1}[2,[1,3]][1,3])\nonumber\\
&&~~~~~~~~ (p_{23}^{ - 1}3p_{21}^{ - 1}p_{23}^{ - 1} (1 - p_{12}p_{21})[1,3] +  p_{21}^{ - 1}p_{23}^{ - 1}p_{12}^{ - 1}p_{13}^{ - 1}p_{21}^{ - 1}p_{23}^{ - 1}
 (1 - p_{12}p_{21})[1,3]3)\nonumber\\
&&~~~~~ +  p_{23}^{ - 7}p_{13}^{ - 7}p_{21}^{ - 2} (1 - \xi)^{5}\xi p_{21}^{ - 1}p_{23}^{ - 1} (1 - p_{12}p_{21}) ([1,3][2,[1,3]]
 +  p_{21}^{ - 1}p_{22}^{ - 1}p_{23}^{ - 1}[2,[1,3]][1,3])3[1,3]\nonumber\\
&&~~~~= - p_{23}^{ - 6}p_{13}^{ - 5} (1 - \xi)^{6}\xi^{2}p_{21}^{ - 2}3[1,3] ([1,3][2,[1,3]] +  p_{21}^{ - 1}p_{22}^{ - 1}p_{23}^{ - 1}[2,[1,3]][1,3])\nonumber\\
&&~~~~~ +  p_{13}^{ - 5}p_{21}^{ - 2}p_{23}^{ - 6}\xi^{2} (1 - \xi)^{6}3[1,3] ([1,3][2,[1,3]] - p_{21}^{ - 1}p_{23}^{ - 1}[2,[1,3]][1,3])\nonumber\\
&&~~~~~ - p_{12}^{ - 1}p_{13}^{ - 7}p_{21}^{ - 4}p_{23}^{ - 8}\xi^{2} (1 - \xi)^{6}
 ([1,3][2,[1,3]] - p_{21}^{ - 1}p_{23}^{ - 1}[2,[1,3]][1,3])3[1,3]\nonumber\\
&&~~~~~ +  p_{23}^{ - 8}p_{13}^{ - 7}p_{21}^{ - 3} (1 - \xi)^{5}\xi  (1 - p_{12}p_{21}) ([1,3][2,[1,3]]
 +  p_{21}^{ - 1}p_{22}^{ - 1}p_{23}^{ - 1}[2,[1,3]][1,3])3[1,3]=0.\nonumber
\end{eqnarray}
Then $[[[1,2],[[1,3],2]],[1,3]]^2=0$.

\vskip.1in
 {\rm (iii)}   Assume generalized Dynkin diagram of braided vector space $V$ is the following condition:

 $\begin{picture} (100,    20) \put (27,   1){\makebox (0,
0)[t]{$\bullet$}} \put (60,    1){\makebox (0,   0)[t]{$\bullet$}}
\put (93,    1){\makebox (0,    0)[t]{$\bullet$}} \put (28,
 - 1){\line (1,
0){33}} \put (61,     - 1){\line (1,    0){30}} \put (18,    7){q} \put (35,    6){$q^{ - 1}$} \put (58, 7){$ - 1$}  \put (75,    6){$r^{ - 1}$}  \put (93,    7){$r$} \put (110,    1){$, $}\ \ \ \put (120,    1){$ q \in R_m,    r \in R_{m'}, q\not= r; m, m' >1$.} \end{picture}$ \\
Then $D (V)=\{ [x_1]; [x_2]; [x_3];  [x_1,x_2];  [x_1,x_3];
[[x_1,x_3],x_2]$; $[[x_1,x_2],[x_1,x_3]]\}$  by  Theorem 3.3(i) in \cite {WZZ14}. Now we show  that every hard super-letter is nilpotent.

 $<y_1,[1,3][1,2]>= (p_{13}^{ - 1} - p_{31})3[1,2] +  p_{11}^{ - 1}p_{13}^{ - 1} (p_{12}^{ - 1} - p_{21})[1,3]2$,

$<y_{2}y_1,[1,3][1,2]>=p_{11}^{ - 1}p_{13}^{ - 1}p_{21}^{ - 1}p_{23}^{ - 1} (p_{12}^{ - 1} - p_{21})[1,3]$,

$<y_{3}y_1,[1,3][1,2]>= (p_{13}^{ - 1} - p_{31})[1,2]$,

$<y_1,[1,2][1,3]>= (p_{12}^{ - 1} - p_{21})2[1,3] +  p_{11}^{ - 1}p_{12}^{ - 1} (p_{13}^{ - 1} - p_{31})[1,2]3$,

$<y_{2}y_1,[1,2][1,3]>= (p_{12}^{ - 1} - p_{21})[1,3]$,

$<y_{3}y_1,[1,2][1,3]>=p_{11}^{ - 1}p_{12}^{ - 1}p_{31}^{ - 1}p_{32}^{ - 1}
 (p_{13}^{ - 1} - p_{31})[1,2]$,

$[[1,2],[1,3]]^{k}= ([1,3][1,2])^{k} +   ( - p_{11}p_{12}p_{31}p_{32})^{k} ([1,2][1,3])^{k}$,
\begin{eqnarray*}
&&<y_1, ([1,3][1,2])^{k}>=\sum \limits _{i=0}^{k - 1} (p_{13}p_{12})^{ - i} ([1,3][1,2])^{i}<y_1,[1,3][1,2]> ([1,3][1,2])^{k - 1 - i}\\
&&~~~~=\sum \limits _{i=0}^{k - 1} (p_{13}p_{12})^{ - i} ([1,3][1,2])^{i} \{ (p_{13}^{ - 1} - p_{31})3[1,2]\\
&&~~~~~ +  p_{11}^{ - 1}p_{13}^{ - 1} (p_{12}^{ - 1} - p_{21})[1,3]2\}
 ([1,3][1,2])^{k - 1 - i},
\end{eqnarray*}
\begin{eqnarray*}
&&<y_{3}y_1, ([1,3][1,2])^{k}>\\
&&~~~~=\sum \limits _{i=0}^{k - 1} (p_{13}p_{12})^{ - i}
(p_{31}^{2}p_{33}p_{32})^{ - i} ([1,3][1,2])^{i} ( (p_{13}^{ - 1} - p_{31})[1,2]) ([1,3][1,2])^{k - 1 - i}\\
&&~~~~= (p_{13}^{ - 1} - p_{31})[1,2] ([1,3][1,2])^{k - 1},
\end{eqnarray*}
\begin{eqnarray*}
&&<y_{2}y_1, ([1,3][1,2])^{k}>\\
&&~~~=\sum \limits _{i=0}^{k - 1} (p_{13}p_{12})^{ - i} (p_{21}^{2}p_{23}p_{22})^{ - i} ([1,3][1,2])^{i}
 (p_{11}^{ - 1}p_{13}^{ - 1}p_{21}^{ - 1}p_{23}^{ - 1}
(p_{12}^{ - 1} - p_{21})[1,3]) ([1,3][1,2])^{k - 1 - i}\\
&&~~~~= (p_{13}p_{12})^{ - n +  1} (p_{21}^{2}p_{23}p_{22})^{ - n +  1} ([1,3][1,2])^{k - 1}
 (p_{11}^{ - 1}p_{13}^{ - 1}p_{21}^{ - 1}p_{23}^{ - 1} (p_{12}^{ - 1} - p_{21})[1,3])\\
&&~~~~= -  (p_{13}p_{21}p_{23})^{ - n} (p_{12}^{ - 1} - p_{21}) ([1,3][1,2])^{k - 1}[1,3],
\end{eqnarray*}
\begin{eqnarray*}
&&<y_1, ([1,2][1,3])^{k}>=\sum \limits _{i=0}^{k - 1} (p_{12}p_{13})^{ - i} ([1,2][1,3])^{i}<y_1,[1,2][1,3]> ([1,2][1,3])^{k - 1 - i}\\
&&~~~=\sum \limits _{i=0}^{k - 1} (p_{12}p_{13})^{ - i} ([1,2][1,3])^{i} ( (p_{12}^{ - 1} - p_{21})2[1,3] +  p_{11}^{ - 1}p_{12}^{ - 1} (p_{13}^{ - 1} - p_{31})[1,2]3)
 ([1,2][1,3])^{k - 1 - i},
\end{eqnarray*}
\begin{eqnarray*}
&&<y_{2}y_1, ([1,2][1,3])^{k}>\\
&&~~~~=\sum \limits _{i=0}^{k - 1} (p_{12}p_{13})^{ - i} (p_{21}^{2}p_{22}p_{23})^{ - i} ([1,2][1,3])^{i}
( (p_{12}^{ - 1} - p_{21})[1,3]) ([1,2][1,3])^{k - 1 - i}\\
&&~~~~= (p_{12}^{ - 1} - p_{21})[1,3] ([1,2][1,3])^{k - 1},
\end{eqnarray*}
\begin{eqnarray*}
&&<y_{3}y_1, ([1,2][1,3])^{k}>\\
&&~~~~=\sum \limits _{i=0}^{k - 1} (p_{12}p_{13})^{ - i} (p_{31}^{2}p_{32}p_{33})^{ - i} ([1,2][1,3])^{i}
 (p_{11}^{ - 1}p_{12}^{ - 1}p_{31}^{ - 1}p_{32}^{ - 1} (p_{13}^{ - 1} - p_{31})[1,2])
([1,2][1,3])^{k - 1 - i}\\
&&~~~~= (p_{12}p_{13})^{ - n +  1} (p_{31}^{2}p_{32}p_{33})^{ - n +  1} ([1,2][1,3])^{k - 1}
 (p_{11}^{ - 1}p_{12}^{ - 1}p_{31}^{ - 1}p_{32}^{ - 1} (p_{13}^{ - 1} - p_{31})[1,2])\\
&&~~~~= -  (p_{12}p_{31}p_{32})^{ - n} (p_{13}^{ - 1} - p_{31}) ([1,2][1,3])^{k - 1}[1,2],
\end{eqnarray*}
\begin{eqnarray*}
&&<y_1,[1,2] ([1,3][1,2])^{k - 1}>= (p_{12}^{ - 1} - p_{21})2 ([1,3][1,2])^{k - 1}
+  p_{11}^{ - 1}p_{12}^{ - 1}[1,2]\\
&&~~~~\sum \limits _{i=0}^{k - 2} (p_{13}p_{12})^{ - i} ([1,3][1,2])^{i}
( (p_{13}^{ - 1} - p_{31})3[1,2] +  p_{11}^{ - 1}p_{13}^{ - 1} (p_{12}^{ - 1} - p_{21})[1,3]2)
 ([1,3][1,2])^{k - 2 - i},
\end{eqnarray*}
\begin{eqnarray*}
&&<y_{2}y_1,[1,2] ([1,3][1,2])^{k - 1}>= (p_{12}^{ - 1} - p_{21}) ([1,3][1,2])^{k - 1} +  p_{11}^{ - 1}p_{12}^{ - 1}[1,2]\\
&&~~~~\sum \limits _{i=0}^{k - 2} (p_{21}^{2}p_{22}p_{23})^{ - i} (p_{13}p_{12})^{ - i}
([1,3][1,2])^{i} (p_{11}^{ - 1}p_{13}^{ - 1}p_{21}^{ - 1}p_{23}^{ - 1} (p_{12}^{ - 1} - p_{21})[1,3])
 ([1,3][1,2])^{k - 2 - i}\\
&&~~~~= (p_{12}^{ - 1} - p_{21}) ([1,3][1,2])^{k - 1} +  p_{11}^{ - 1}p_{12}^{ - 1}[1,2]
 (p_{21}^{2}p_{22}p_{23})^{ - n +  2} (p_{13}p_{12})^{ - n +  2} \\
&&~~~~~~~([1,3][1,2])^{k - 2}p_{11}^{ - 1}p_{13}^{ - 1}p_{21}^{ - 1}p_{23}^{ - 1}
(p_{12}^{ - 1} - p_{21})[1,3]\\
&&~~~~= (p_{12}^{ - 1} - p_{21}) ([1,3][1,2])^{k - 1} +    (p_{21}p_{13}p_{23})^{ - n +  1}p_{12}^{ - 1} (p_{12}^{ - 1} - p_{21}) ([1,2][1,3])^{k - 1},\\
&&<y_{3}y_1,[1,2] ([1,3][1,2])^{k - 1}>=0,
\end{eqnarray*}
\begin{eqnarray*}
&&<y_1,[1,3] ([1,2][1,3])^{k - 1}>= (p_{13}^{ - 1} - p_{31})3 ([1,2][1,3])^{k - 1} +  p_{11}^{ - 1}p_{13}^{ - 1}[1,3]\\
&&~~~~\sum \limits _{i=0}^{k - 2} (p_{12}p_{13})^{ - i} ([1,2][1,3])^{i}
( (p_{12}^{ - 1} - p_{21})2[1,3] +  p_{11}^{ - 1}p_{12}^{ - 1} (p_{13}^{ - 1} - p_{31})[1,2]3)
 ([1,2][1,3])^{k - 2 - i},
\end{eqnarray*}
\begin{eqnarray*}
&&<y_{3}y_1,[1,3] ([1,2][1,3])^{k - 1}>= (p_{13}^{ - 1} - p_{31}) ([1,2][1,3])^{k - 1} +  p_{11}^{ - 1}p_{13}^{ - 1}[1,3]\\
&&~~~~\sum \limits _{i=0}^{k - 2} (p_{31}^{2}p_{33}p_{32})^{ - i} (p_{12}p_{13})^{ - i}
([1,2][1,3])^{i} (p_{11}^{ - 1}p_{12}^{ - 1}p_{31}^{ - 1}p_{32}^{ - 1} (p_{13}^{ - 1} - p_{31})[1,2])
 ([1,2][1,3])^{k - 2 - i}\\
&&~~~~= (p_{13}^{ - 1} - p_{31}) ([1,2][1,3])^{k - 1}\\
&&~~~~~ +  p_{11}^{ - 1}p_{13}^{ - 1}[1,3]
 (p_{31}^{2}p_{33}p_{32})^{ - n +  2} (p_{12}p_{13})^{ - n +  2} ([1,2][1,3])^{k - 2}p_{11}^{ - 1}p_{12}^{ - 1}p_{31}^{ - 1}p_{32}^{ - 1}(p_{13}^{ - 1} - p_{31})[1,2]\\
&&~~~~= (p_{13}^{ - 1} - p_{31}) ([1,2][1,3])^{k - 1} +
 (p_{31}p_{12}p_{32})^{ - n +  1}p_{13}^{ - 1} (p_{13}^{ - 1} - p_{31}) ([1,3][1,2])^{k - 1},
\end{eqnarray*}
\begin{eqnarray*}
&&<y_{2}y_1,[1,3] ([1,2][1,3])^{k - 1}>=0,\\
&&<y_{2}y_{1}y_{2}y_1, ([1,2][1,3])^{k}>=0,\\
&&<y_{2}y_{1}y_{2}y_1, ([1,3][1,2])^{k}>=0,\\
&&<y_{3}y_{1}y_{3}y_1, ([1,2][1,3])^{k}>=0,\\
&&<y_{3}y_{1}y_{3}y_1, ([1,3][1,2])^{k}>=0,\\
&&<y_{2}y_{1}y_{3}y_1, ([1,2][1,3])^{k}>=<y_{2}y_{1}, -  (p_{12}p_{31}p_{32})^{ - n} (p_{13}^{ - 1} - p_{31}) ([1,2][1,3])^{k - 1}[1,2]>\\
&&~~~~= -  (p_{12}p_{31}p_{32})^{ - n} (p_{13}^{ - 1} - p_{31})<y_{2}y_{1}, ([1,2][1,3])^{k - 1}[1,2]>,
\end{eqnarray*}
\begin{eqnarray*}
&&<y_{2}y_{1}y_{3}y_1, ([1,3][1,2])^{k}>=<y_{2}y_{1}, (p_{13}^{ - 1} - p_{31})[1,2] ([1,3][1,2])^{k - 1}>\\
&&~~~~= (p_{13}^{ - 1} - p_{31})<y_{2}y_{1},[1,2] ([1,3][1,2])^{k - 1}>,
\end{eqnarray*}
\begin{eqnarray*}
&&<y_{3}y_{1}y_{2}y_1, ([1,3][1,2])^{k}>=<y_{3}y_{1}, -  (p_{13}p_{21}p_{23})^{ - n} (p_{12}^{ - 1} - p_{21}) ([1,3][1,2])^{k - 1}[1,3]>\\
&&~~~~= -  (p_{13}p_{21}p_{23})^{ - n} (p_{12}^{ - 1} - p_{21})<y_{3}y_{1}, ([1,3][1,2])^{k - 1}[1,3]>,
\end{eqnarray*}
\begin{eqnarray*}
&&<y_{3}y_{1}y_{2}y_1, ([1,2][1,3])^{k}>=<y_{3}y_{1}, (p_{12}^{ - 1} - p_{21})[1,3] ([1,2][1,3])^{k - 1}>\\
&&~~~~= (p_{12}^{ - 1} - p_{21})<y_{3}y_{1},[1,3] ([1,2][1,3])^{k - 1}>,
\end{eqnarray*}
\begin{eqnarray*}
&&<y_{2}y_{1}y_{3}y_1,[[1,2],[1,3]]^{k}>=<y_{2}y_{1}y_{3}y_1, ([1,3][1,2])^{k} +   ( - p_{11}p_{12}p_{31}p_{32})^{k} ([1,2][1,3])^{k}>\\
&&~~~~=<y_{2}y_{1}y_{3}y_1, ([1,3][1,2])^{k}> +  (p_{12}p_{31}p_{32})^{k}<y_{2}y_{1}y_{3}y_1, ([1,2][1,3])^{k}>\\
&&~~~~= (p_{13}^{ - 1} - p_{31})<y_{2}y_{1},[1,2] ([1,3][1,2])^{k - 1}> \\
&&~~~~~- (p_{12}p_{31}p_{32})^{k}
 (p_{12}p_{31}p_{32})^{ - n} (p_{13}^{ - 1} - p_{31})<y_{2}y_{1}, ([1,2][1,3])^{k - 1}[1,2]>=0,
\end{eqnarray*}
\begin{eqnarray*}
&&<y_{3}y_{1}y_{2}y_1,[[1,2],[1,3]]^{k}>=<y_{3}y_{1}y_{2}y_1, ([1,3][1,2])^{k} +   ( - p_{11}p_{12}p_{31}p_{32})^{k} ([1,2][1,3])^{k}>\\
&&~~~~=<y_{3}y_{1}y_{2}y_1, ([1,3][1,2])^{k}> +  (p_{12}p_{31}p_{32})^{k}<y_{3}y_{1}y_{2}y_1, ([1,2][1,3])^{k}>\\
&&~~~~= -  (p_{13}p_{21}p_{23})^{ - n} (p_{12}^{ - 1} - p_{21})<y_{3}y_{1}, ([1,3][1,2])^{k - 1}[1,3]> \\
&&~~~~~+ (p_{12}p_{31}p_{32})^{k} (p_{12}^{ - 1} - p_{21})<y_{3}y_{1},[1,3] ([1,2][1,3])^{k - 1}>\\
&&~~~~= (p_{13}p_{21}p_{23})^{ - n} ( (p_{22}p_{33})^{ - n} - 1) (p_{12}^{ - 1} - p_{21})
<y_{3}y_{1}, ([1,3][1,2])^{k - 1}[1,3]>.
\end{eqnarray*}
Then $[[1,2],[1,3]]^{k}=0,~~~~  k=\frac {m^2 m'{}^2}{ (m,   m') }$. \hfill $\Box$

\begin {Lemma} \label {8'} (Theorem 3.1(1) in \cite {AS00}, and \cite {He05}) There does not exist any $m$-infinity element in  $\mathfrak B (V)$
when $\Delta  (\mathfrak B (V))$ is an arithmetic root system and $V$ is of finite Cartan type with
${\rm ord }  (q_{ij}) < \infty$ for $1\le i, j \le n$ or $\dim V =n <3.$ \end {Lemma}

\subsection {The hard super-letters of $\mathbb B(V)$ with $\dim V =2$}
We give  $D  := \{u_1, u_2, \cdots, u_r\}$ for all hard super-letters of connected braided vector space $V$  of diagonal type with rank $2$ and height $h_i$ for hard super-letter $u_i$ using Definition 2 and Appendix A in \cite {He05}.
Let $a = 112, b= 122, c= 1112, d = 11212, e = 11112, f = 1112112, g = 1121212, h = 111112,  i = 111121112, j = 11121112112$, $k = 1112112112,l =11211212, m = 112121121212, n =112121212, p = 11211211212, q = 1121121211212, r = 12 $. $h_u = {\rm ord }  (p _{u u})$ for any $u\in D (V)$.

$\hbox {T{2}}.$  $D= \{1,2,r\}$.

$\hbox {T{3}}.$  $D=\{1,2,r,a\}$.

$\hbox {T{4}}.$  $D=\{1,2,r,b,a\}$.

$\hbox {T{5}}.$  $D=\{1,2,r,a,d\}$.

$\hbox {T{6}}.$  $D=\{1,2,r,b,a,d\}$.

$\hbox {T{7}}.$  $D=\{1,2,r,a,c\}$.

$\hbox {T{8}}.$  $D=\{1,2,r,a,d,c\}$.

$\hbox {T{9}}.$  $D=\{1,2,r,a,d,g\}$.

$\hbox {T{10}}.$   $D=\{1,2,r,b,a,d,g,c\}$.

$\hbox {T{11}}.$  $D=\{1,2,r,a,d,g,l,c\}$.

$\hbox {T{12}}.$  $D=\{1,2,r,b,a,d,c,f\}$.

$\hbox {T{13}}.$  $D=\{1,2,r,a,d,l,c,f\}$.

$\hbox {T{14}}.$   $D=\{1,2,r,a,c,e\}$.

$\hbox {T{15}}.$  $D=\{1,2,r,a,d,g,c,e\}$.

$\hbox {T{16}}.$  $D=\{1,2,r,a,d,c,f,e\}$.

$\hbox {T{17}}.$  $D=\{1,2,r,a,d,g,n,l\}$.

$\hbox {T{18}}.$  $D=\{1,2,r,a,d,g,l,q\}$

$\hbox {T{19}}.$  $D=\{1,2,r,a,d,g,n,m,l,q,p,c\}$.

$\hbox {T{20}}$  $D=\{1,2,r,a,c,f,j,e\}$.

$\hbox {T{21}}$  $D=\{1,2,r,a,c,f,e,h\}$.

$\hbox {T{22}}$  $D=\{1,2,r,a,d,c,f,k,j,e,i,h\}$.

 (Row $2$)$\hbox {T{2}}$  Table A.1 in \cite  {He05} : we set $u_{1}=x_{2},u_{2}=[r],u_{3}=x_{1}$.
 $\{p_{u_i, u_i} \mid 1\le i \le 3\} = \{q, q, q \}$.

 (Row $3$) $\hbox {T{2} (1)}$  Table A.1 in \cite  {He05} : we set $u_{1}=x_{2},u_{2}=[r],u_{3}=x_{1}$. $\{p_{u_i, u_i} \mid 1\le i \le 3\} = \{ -1, -1, q \}$.

$\hbox {T{2} (2)}$  Table A.1 in \cite  {He05} : we set $u_{1}=x_{2},u_{2}=[r], u_{3}=x_{1}$.
 $\{p_{u_i, u_i} \mid 1\le i \le 3\} = \{-1, q, -1 \}$.

 (Row $4$) $\hbox {T{3}}$  Table A.1 in \cite  {He05} : we set $u_{1}=x_{2},u_{2}=[r], u_{3}=[a],
u_{4}=x_{1}$.  $\{p_{u_i, u_i} \mid 1\le i \le 4\} = \{q^2, q, q^2,  q \}$.

 (Row $5$) $\hbox {T{3} (1)}$  Table A.1 in \cite  {He05} : we set $u_{1}=x_{2},u_{2}=[r],u_{3}=[a],
u_{4}=x_{1}$.  $\{p_{u_i, u_i} \mid 1\le i \le 4\} = \{-1, -q^{-1}, -1, q \}$.

$\hbox {T{3} (2)}$  Table A.1 in \cite  {He05}: we set $u_{1}=x_{2},u_{2}=[r], u_{3}=[a],
u_{4}=x_{1}$.  $\{p_{u_i, u_i} \mid 1\le i \le 4\} = \{-1, q, -1, - q^{-1} \}$.

 (Row $6$) $\hbox {T{3} (1)}$ Table A.1 in \cite  {He05} : we set $u_{1}=x_{2},u_{2}=[r], u_{3}=[a],
u_{4}=x_{1}$.  $\{p_{u_i, u_i} \mid 1\le i \le 4\} = \{q, \xi, \xi q^{-1}, \xi \}$.

$\hbox {T{3} (2)}$  Table A.1 in \cite  {He05}: we set $u_{1}=x_{2},u_{2}=[r], u_{3}=[a],
u_{4}=x_{1}$.  $\{p_{u_i, u_i} \mid 1\le i \le 4\} = \{ \xi q^{-1}, \xi, q, \xi \}$.

 (Row $7$) $\hbox {T{3} (1)}$  Table A.1 in \cite  {He05} : we set $u_{1}=x_{2},u_{2}=[r],u_{3}=[a],
u_{4}=x_{1}$.  $\{p_{u_i, u_i} \mid 1\le i \le 4\} = \{-1,\xi^{-1}, -1, \xi \}$.

$\hbox {T{3} (2)}$  Table A.1 in \cite  {He05} : we set $u_{1}=x_{2},u_{2}=[r],u_{3}=[a],
u_{4}=x_{1}$.  $\{p_{u_i, u_i} \mid 1\le i \le 4\} = \{-1, \xi, -1, \xi^{-1} \}$.

 (Row $8$) $\hbox {T{4}}$  Table A.1 in \cite  {He05} : we set $u_{1}=x_{2},u_{2}=[b],u_{3}=[r],
u_{4}=[a],
u_{5}=x_{1}$.  $\{p_{u_i, u_i} \mid 1\le i \le 5\} = \{ -\xi^{2}, -1, -\xi^3, -1, -\xi^{-2} \}$.

$\hbox {T{5} (1)}$  Table A.1 in \cite  {He05} : we set $u_{1}=x_{2},u_{2}=[r],u_{3}=[d],
u_{4}=[a],u_{5}=x_{1}$.  $\{p_{u_i, u_i} \mid 1\le i \le 5\} = \{-1, -\xi^3, -1, -\xi^2, -\xi^{-2} \}$.

$\hbox {T{5} (2)}$  Table A.1 in \cite  {He05} : we set $u_{1}=x_{2},u_{2}=[r],u_{3}=[d],
u_{4}=[a],u_{5}=x_{1}$.  $\{p_{u_i, u_i} \mid 1\le i \le 5\} = \{-1, -\xi^3, -1, -\xi^{-2}, -\xi^2\}$.

$\hbox {T{7} (1)}$  Table A.1 in \cite  {He05} : we set $u_{1}=x_{2},u_{2}=[r],u_{3}=[a],
u_{4}=[c],u_{5}=x_{1}$.  $\{p_{u_i, u_i} \mid 1\le i \le 5\} = \{-1, -\xi ^{-2}, -\xi^2, -1, -\xi ^3 \}$.

$\hbox {T{7} (2)}$  Table A.1 in \cite  {He05} : we set $u_{1}=x_{2},u_{2}=[r],u_{3}=[a],
u_{4}=[c],u_{5}=x_{1}$.  $\{p_{u_i, u_i} \mid 1\le i \le 5\} = \{-1, -\xi^2, -\xi ^{-2}, -1-\xi^3\}$.

 (Row $9$) $\hbox {T{4}}$  Table A.1 in \cite  {He05} : we set $u_{1}=x_{2},u_{2}=[b],u_{3}=[r],
u_{4}=[a],
u_{5}=x_{1}$.  $\{p_{u_i, u_i} \mid 1\le i \le 5\} = \{-\xi ^2, -1, -\xi^{-1}, -1, -\xi ^2 \}$.

$\hbox {T{5}}$  Table A.1 in \cite  {He05} : we set $u_{1}=x_{2},u_{2}=[r],u_{3}=[d],
u_{4}=[a],u_{5}=x_{1}$.  $\{p_{u_i, u_i} \mid 1\le i \le 5\} = \{-1, -\xi^{-1}, -1, -\xi ^ 2, -\xi ^2 \}$.

$\hbox {T{7}}$  Table A.1 in \cite  {He05} : we set $u_{1}=x_{2},u_{2}=[r],u_{3}=[a],
u_{4}=[c],u_{5}=x_{1}$.  $\{p_{u_i, u_i} \mid 1\le i \le 5\} = \{-1, -\xi^2, -\xi ^2, -1, -\xi ^{-1} \}$.

 (Row $10$) $\hbox {T{6}}$  Table A.1 in \cite  {He05} : we set $u_{1}=x_{2},u_{2}=[b],u_{3}=[r],
u_{4}=[d],u_{5}=[a],
u_{6}=x_{1}$.  $\{p_{u_i, u_i} \mid 1\le i \le 6\} = \{\xi ^3, -1, -\xi ^2, -1, \xi ^3 -\xi\}$.

$\hbox {T{9}}$  Table A.1 in \cite  {He05} : we set $u_{1}=x_{2},u_{2}=[r],u_{3}=[g],
u_{4}=[d],u_{5}=[a],
u_{6}=x_{1}$.  $\{p_{u_i, u_i} \mid 1\le i \le 6\} = \{-1, -\xi ^2, -1, \xi ^3, -\xi, \xi ^3 \}$.

$\hbox {T{14}}$  Table A.1 in \cite  {He05} : we set $u_{1}=x_{2},u_{2}=[r],u_{3}=[a],
u_{4}=[c],u_{5}=[e],
u_{6}=x_{1}$.   $\{p_{u_i, u_i} \mid 1\le i \le 6\} = \{-1, \xi ^3, -\xi, \xi^3, -1, -\xi^2 \}$.

 (Row $11$) $\hbox {T{8}}$  Table A.1 in \cite  {He05} : we set $u_{1}=x_{2},u_{2}=[r],u_{3}=[d],
u_{4}=[a],u_{5}=[c],
u_{6}=x_{1}$.  $\{p_{u_i, u_i} \mid 1\le i \le 6\} = \{q^3, q, q^3, q , q^3, q \}$.

 (Row $12$) $\hbox {T{8} (1)}$  Table A.1 in \cite  {He05} : we set $u_{1}=x_{2},u_{2}=[r],u_{3}=[d],
u_{4}=[a],u_{5}=[c],
u_{6}=x_{1}$.  $\{p_{u_i, u_i} \mid 1\le i \le 6\} = \{\xi ^{-1}, \xi ^2, -1, \xi, -1, \xi ^2 \}$.

$\hbox {T{8} (2)}$  Table A.1 in \cite  {He05} : we set $u_{1}=x_{2},u_{2}=[r],u_{3}=[d],
u_{4}=[a],u_{5}=[c],
u_{6}=x_{1}$.   $\{p_{u_i, u_i} \mid 1\le i \le 6\} = \{-1, \xi, -1, \xi^2, \xi^{-1}, \xi ^2 \}$.

$\hbox {T{8} (3)}$  Table A.1 in \cite  {He05} : we set $u_{1}=x_{2},u_{2}=[r],u_{3}=[d],
u_{4}=[a],u_{5}=[c],
u_{6}=x_{1}$. $\{p_{u_i, u_i} \mid 1\le i \le 6\} = \{-1, \xi^{2}, \xi^{-1}, \xi^2, -1, \xi \}$.

 (Row $13$) $\hbox {T{10}}$  Table A.1 in \cite  {He05} : we set $u_{1}=x_{2},u_{2}=[b],u_{3}=[r],
u_{4}=[g],u_{5}=[d],u_{6}=[a],u_{7}=[c],
u_{8}=x_{1}$.  $\{p_{u_i, u_i} \mid 1\le i \le 8\} = \{-\xi^{-4}, -1, \xi, -1, -\xi ^{-4}, \xi ^6, \xi^{-1}, \xi ^6 \}$.

$\hbox {T{13}}$ Table A.1 in \cite  {He05} : we set $u_{1}=x_{2},u_{2}=[r],u_{3}=[d],
u_{4}=[l],u_{5}=[a],u_{6}=[f],u_{7}=[c],
u_{8}=x_{1}$.  $\{p_{u_i, u_i} \mid 1\le i \le 8\} = \{\xi^{-1}, \xi ^6, -\xi ^{-4}, -1, \xi, -1, -\xi^{-4}, \xi^6 \}$.

$\hbox {T{17}}$  Table A.1 in \cite  {He05}: we set $u_{1}=x_{2},u_{2}=[r],u_{3}=[n],
u_{4}=[g],u_{5}=[d],u_{6}=[l],u_{7}=[a],
u_{8}=x_{1}$.  $\{p_{u_i, u_i} \mid 1\le i \le 8\} = \{-1, \xi, -1, -\xi^{-4}, \xi ^6, \xi ^{-1}, \xi ^6, -\xi ^{-4} \}$.

$\hbox {T{21}}$ Table A.1 in \cite  {He05} : we set $u_{1}=x_{2},u_{2}=[r],u_{3}=[a],
u_{4}=[f],u_{5}=[c],u_{6}=[e],u_{7}=[h],
u_{8}=x_{1}$.  $\{p_{u_i, u_i} \mid 1\le i \le 8\} = \{-1, -\xi ^{-4}, \xi ^6, \xi ^{-1}, \xi ^6, -\xi ^{-4}, -1, \xi \}$.

 (Row $14$) $\hbox {T{11}}$  Table A.1 in \cite  {He05} : we set $u_{1}=x_{2},u_{2}=[r],u_{3}=[g],
u_{4}=[d],u_{5}=[l],u_{6}=[a],u_{7}=[c],
u_{8}=x_{1}$.   $\{p_{u_i, u_i} \mid 1\le i \le 8\} = \{-1, -\xi ^{-2}, -1, \xi, -1, -\xi^{-2}, -1, \xi \}$.

$\hbox {T{16}}$ Table A.1 in \cite  {He05} : we set $u_{1}=x_{2},u_{2}=[r],u_{3}=[d],
u_{4}=[a],u_{5}=[f],u_{6}=[c],u_{7}=[e],
u_{8}=x_{1}$.  $\{p_{u_i, u_i} \mid 1\le i \le 8\} = \{-1, \xi, -1, -\xi^{-2}, -1, \xi, -1, -\xi ^{-2} \}$.

 (Row $15$) $\hbox {T{11} (1)}$  Table A.1 in \cite  {He05} : we set $u_{1}=x_{2},u_{2}=[r],u_{3}=[g],
u_{4}=[d],u_{5}=[l],u_{6}=[a],u_{7}=[c],
u_{8}=x_{1}$.   $\{p_{u_i, u_i} \mid 1\le i \le 8\} = \{-1, -\xi ^{-2}, -1, -\xi, -1, -\xi^{-2}, -1, \xi \}$.

$\hbox {T{11} (2)}$  Table A.1 in \cite  {He05} : we set $u_{1}=x_{2},u_{2}=[r],u_{3}=[g],
u_{4}=[d],u_{5}=[l],u_{6}=[a],u_{7}=[c],
u_{8}=x_{1}$.
 $\{p_{u_i, u_i} \mid 1\le i \le 8\} = \{ -1, -\xi ^{-2}, -1, \xi, -1, -\xi^{-2}, -1, -\xi\}$.

$\hbox {T{16} (1)}$  Table A.1 in \cite  {He05} : we set $u_{1}=x_{2},u_{2}=[r],u_{3}=[d],
u_{4}=[a],u_{5}=[f],u_{6}=[c],u_{7}=[e],
u_{8}=x_{1}$.   $\{p_{u_i, u_i} \mid 1\le i \le 8\} = \{-1, \xi, -1, -\xi ^{-2}, -1, -\xi, -1, -\xi ^{-2}\}$.

$\hbox {T{16} (2)}$  Table A.1 in \cite  {He05} : we set $u_{1}=x_{2},u_{2}=[r],u_{3}=[d],
u_{4}=[a],u_{5}=[f],u_{6}=[c],u_{7}=[e],
u_{8}=x_{1}$.  $\{p_{u_i, u_i} \mid 1\le i \le 8\} = \{-1, -\xi, -1, -\xi^{-2}, -1, \xi, -1, -\xi^{-2} \}$.

 (Row $16$) $\hbox {T{12}}$  Table A.1 in \cite  {He05} : we set $u_{1}=x_{2},u_{2}=[b],u_{3}=[r],
u_{4}=[d],u_{5}=[a],u_{6}=[f],u_{7}=[c],
u_{8}=x_{1}$.  $\{p_{u_i, u_i} \mid 1\le i \le 8\} = \{\xi ^5, -1, \xi^3, -\xi^{-4}, \xi ^3, -1, \xi ^5, -\xi \}$.

$\hbox {T{15}}$  Table A.1 in \cite  {He05} : we set $u_{1}=x_{2},u_{2}=[r],u_{3}=[g],
u_{4}=[d],u_{5}=[a],u_{6}=[c],u_{7}=[e],
u_{8}=x_{1}$.  $\{p_{u_i, u_i} \mid 1\le i \le 8\} = \{ -\xi^{-4}, \xi ^3, -1, \xi ^5,
-\xi, \xi ^5, -\xi, \xi ^5, -1, \xi ^3 \}$.

$\hbox {T{18}}$  Table A.1 in \cite  {He05} : we set $u_{1}=x_{2},u_{2}=[r],u_{3}=[g],
u_{4}=[d],u_{5}=[q],u_{6}=[l],u_{7}=[a],
u_{8}=x_{1}$.  $\{p_{u_i, u_i} \mid 1\le i \le 8\} = \{-1, \xi ^3, -\xi ^{-4}, \xi ^3, -1, \xi ^5, -\xi,  \xi ^5 \}$.

$\hbox {T{20}}$ Table A.1 in \cite  {He05} : we set $u_{1}=x_{2},u_{2}=[r],u_{3}=[a],
u_{4}=[f],u_{5}=[j],u_{6}=[c],u_{7}=[e],
u_{8}=x_{1}$.  $\{p_{u_i, u_i} \mid 1\le i \le 8\} = \{-1, \xi ^5, -\xi, \xi ^5, -1, \xi ^3, -\xi ^{-4}, \xi ^3 \}$.

 (Row $17$) $\hbox {T{19}}$  Table A.1 in \cite  {He05} : we set $u_{1}=x_{2},u_{2}=[r],u_{3}=[n],u_{4}=[g],
u_{5}=[m],u_{6}=[d],u_{7}=[q],
u_{8}=[l],u_{9}=[p],u_{10}=[a],u_{11}=[c],
u_{12}=x_{1}$.  $\{p_{u_i, u_i} \mid 1\le i \le 12\} = \{-1, -\xi ^{-2}, -1, -\xi, -1, -\xi ^{-2}, -1, -\xi, -1, -\xi ^{-2}, -1, -\xi \}$.

$\hbox {T{22}}$ Table A.1 in \cite  {He05}: we set $u_{1}=x_{2},u_{2}=[r],u_{3}=[d],u_{4}=[a],
u_{5}=[k],u_{6}=[f],u_{7}=[j],
u_{8}=[c],u_{9}=[i],u_{10}=[e],u_{11}=[h],
u_{12}=x_{1}$.  $\{p_{u_i, u_i} \mid 1\le i \le 12\} = \{-1, -\xi, -1, -\xi ^ {-2}, -1, -\xi, -1, -\xi ^{-2}, -1, -\xi, -1, -\xi^{-2} \}$.

\section*{Acknowledgment}
We would like to thank the referee for many constructive suggestions which bring essential improvements to section 2 of the paper. YZZ was partially supported by the Australian Research Council through Discovery-Projects grant DP140101492.

\end {document}